yang@math.duke.edu

\documentstyle[amscd]{amsart}

\def\Z{{\Bbb Z}}
\def\Q{{\Bbb Q}}
\def\R{{\Bbb R}}
\def\C{{\Bbb C}}
\def\P{{\Bbb P}}

\def\calC{{\cal C}}
\def\calD{{\cal D}}

\def\O{{\cal O}}			
\def\H{{\cal H}}
\def\F{{\cal F}}
\def\g{{\frak g}}
\def\h{{\frak h}}
\def\Atilde{\widetilde{A}}
\def\phihat{\hat{\phi}}

\def\dot{\bullet}
\def\vol{\mbox{vol}}

\def\fbar{\overline{f}}
\def\lbar{\overline{l}}

\def\Hcts{H_{\text{cts}}}	
\def\Hp{H_{\text{loc-}L_p}}    	
\def\Hd{H_\calD}         	
\def\Hf{H_\F}

\def\supdot{^\bullet}		
\def\subdot{_\bullet}		

\def\otilde{\widetilde{\O}}
\def\Xbar{{\overline{X}}}
\def\Ytilde{\widetilde{Y}}

\def\bglm{B_\bullet GL_m(\C)}

\def\bgm{B_\bullet G^m}
\def\bum{B_\bullet U^m}
\def\egm{E_\bullet G^m}
\def\eum{E_\bullet U^m}

\def\gen{{\text{gen}}}
\def\Btilde{\widetilde{B}}
\def\Etilde{\widetilde{E}}
\def\cts{{\text{cts}}}

\def\bw{{\cal BW}}
\def\obw{\Omega_{\cal BW}}
\def\dbcomplex{\R(p)_{\cal D}}

\def\blank{{\phantom{X}}}

\newcommand\sgn{\operatorname{sgn}}
\newcommand\spec{\operatorname{Spec}}
\newcommand\im{\operatorname{im}}		
\newcommand\coker{\operatorname{coker}}
\newcommand\cone{\operatorname{Cone}}
\newcommand\Hom{\operatorname{Hom}}
\newcommand\Ext{\operatorname{Ext}}
\newcommand\Alt{\operatorname{Alt}}		
\newcommand\Arg{\operatorname{Arg}}		
\newcommand\Alb{\operatorname{Alb}}		
\renewcommand\Re{\operatorname{Re}}		
\renewcommand\Im{\operatorname{Im}}		

\newtheorem{theorem}{Theorem}[section]
\newtheorem{lemma}[theorem]{Lemma}
\newtheorem{proposition}[theorem]{Proposition}
\newtheorem{corollary}[theorem]{Corollary}

\theoremstyle{definition}
\newtheorem{definition}[theorem]{Definition}
\newtheorem{example}[theorem]{Example}

\theoremstyle{remark}
\newtheorem{remark}[theorem]{Remark}
\newtheorem{variant}[theorem]{Variant}

\begin{document}

\title{Real Grassmann Polylogarithms and Chern Classes}

\author{Richard M.~Hain}

\author{Jun Yang}

\address{Department of Mathematics\\
Duke University\\ Durham, NC 27708-0320}
\email{hain@@math.duke.edu, yang@@math.duke.edu}

\thanks{Both authors were partially supported by grants from the
National
Science Foundation.}

\date{July, 1994}

\maketitle

\section{Introduction}\label{intro}

In this paper we define and prove the existence of real Grassmann
polylogarithms which are the real single-valued analogues of the
Grassmann polylogarithms defined in \cite{hain-macp} and constructed
in \cite{hain-macp,hana-macp,hana-macp_2,hain:generic}. We prove
that if $\eta_X$ is the generic point of a complex algebraic variety
$X$, then the $m$th such polylogarithm represents
the Beilinson Chern class
$$
c_m^B : r_m K_n(\eta_X) \to \Hd^{2m-m}(\eta_X,\R(m)),
$$
where $r_m K_n$ denotes the part of the algebraic $K$-theory coming
from $GL_m$. In particular, we show that there is a {Bloch-Wigner
function} $D_m$ defined on a Zariski open subset of the
grassmannian of $m$ planes in $\C^{2m}$ which satisfies a canonical
$(2m+1)$-term functional equation.

In order to describe our results in more detail, we recall
the definition of Bloch-Wigner functions from \cite[\S
11]{hain-macp}.
Denote the algebra of global holomorphic forms with logarithmic
singularities at infinity on a smooth complex algebraic variety
$X$ by $\Omega^\dot(X)$. A {\it Bloch-Wigner function} on $X$ is
simply a real single-valued function on $X$ that is a polynomial
in real and imaginary parts of multivalued functions on $X$ of the
form
$$
x \mapsto \sum_I\int_{x_0}^x w_{i_1}\dots w_{i_r}
$$
where each $w_{i_j}$ is in $\Omega^1(X)$. For example,
the logarithm and classical dilogarithm $\ln_2(x)$ can be expressed
as iterated integrals,
$$
\log x = \int_1^x \frac{dz}{z},\quad
\ln_2(x) = \int_0^x \frac{dz}{1-z}\frac{dz}{z},
$$
so the single valued logarithm
$$
D_1(x) = \log|x|
$$
and the Bloch-Wigner function
$$
D_2(x) = \Im \ln_2(x) + \log|x|\Arg(1-x)
$$
are both Bloch-Wigner functions. More generally, Ramakrishnan's
single-valued cousins of the classical polylogarithms
$\ln_m(x)$ are Bloch-Wigner functions. The set of Bloch-Wigner
functions on $X$ is an $\R$ algebra, which will be denoted by
$\bw(X)$.

Define the irregularity $q(X)$ of $X$ to be half the first Betti
number of any smooth compactification of $X$.  When $q(X)=0$, the
elements of $\bw(X)$ have a more intrinsic description --- they are
precisely the functions that occur as the matrix entries of period
maps of real variations of mixed Hodge structure over $X$, all of
whose weight graded quotients are constant variations of type
$(p,p)$.

The complex $\Omega_\bw^\dot(X)$ of Bloch-Wigner forms on $X$ is
defined to be the subcomplex of the real de~Rham complex of $X$
that is generated by $\bw(X)$ and the real and imaginary parts
of elements of $\Omega^\dot(X)$ (cf.\ \cite{yang:thesis}.) It is
closed under exterior
differentiation, and is therefore a differential graded algebra.
It has a natural weight filtration $W_\dot$ which comes from the
usual
weight filtrations on $\Omega^\dot(X)$ and on iterated integrals.
It is easy to describe the weight filtration on
$\bw(X)$ when $q(X)=0$: the weight of an iterated integral of length
$l$ is simply $2l$. This filtration induces one on $\bw(X)$
in the natural way, so that
$$
D_1 \in W_2\bw(\C^\ast)\quad\text{and}\quad D_2 \in
W_4\bw(\C-\{0,1\}),
$$
for example.

The Bloch-Wigner (or BW-)cohomology of $X$ is the analogue of the
Deligne cohomology of $X$ constructed using Bloch-Wigner forms
in place of usual forms. Specifically, the Bloch-Wigner cohomology
of $X$ is the cohomology of the complex
$$
\R(m)_\bw\supdot(X) := \cone[F^pW_{2m}\Omega\supdot(X)
\rightarrow W_{2m}\obw\supdot(X)\otimes \R(m-1)][-1],
$$
where the map $\Omega^\dot(X) \to \obw^\dot(X)\otimes\R(m-1)$ is
defined
by taking a complex valued form to its reduction mod $\R(m)$.%
\footnote{Recall that $\R(m)$ is the subgroup $(2\pi i)^m\R$ of $\C$.
There is a standard identification $\C/\R(m)\cong \R(m-1)$.}
The inclusion of $\Omega_\bw^\dot(X)$ into the real log complex
associated to $X$ induces a natural map
$$
H\supdot_\bw(X,\R(m))\to \Hd\supdot(X,\R(m)).
$$

Grassmann polylogarithms are specific elements of the BW-cohomology
of certain simplicial varieties $G^m_\dot$. Recall from
\cite{hain-macp} that the variety $G^m_n$ is defined to be the
Zariski open subset of the grassmannian of $n$ dimensional linear
subspaces of $\P^{m+n}$ which consists of those $n$ planes that
are transverse to each stratum of the union of the coordinate
hyperplanes. Intersecting elements of $G^m_n$ with the $j$th
coordinate hyperplane defines ``face maps''
$$
A_j : G^m_n \to G^m_{n-1}, \qquad j=0,\dots,m+n
$$
These satisfy the usual simplicial identities.
It is natural to place $G^m_n$ in dimension $m+n$ as there are
$m+n+1$ face maps emanating from it. The collection
of the $G^m_n$ with $0\le n \le m$ and the face maps $A_j$ will be
denoted by $G^m_\dot$. It is a truncated simplicial variety.

We can apply the functor $\R(m)_\bw\supdot(\blank)$ to a simplicial
variety $X_\dot$ to obtain a double complex
$\R(m)_\bw\supdot(X_\dot)$.
The homology of the associated total complex will be denoted by
$H_\bw^\dot(X_\dot,\R(m))$.

There is a natural map
$$
\delta : H^{2m}_\bw(G^m\subdot,\R(m)) \to \Omega^m(G^m_0).
$$
A {\it real Grassmann $m$-logarithm} is an element $L_m$ of
$H^{2m}_\bw(G^m\subdot,\R(m))$ whose image under $\delta$ is
the ``volume form''
$$
\vol_m = \frac{dx_1}{x_1}\wedge \dots \wedge \frac{dx_m}{x_m}
$$
on $G^m_0 \cong (\C^\ast)^m$. This is the real analogue of the
definition of a generalized $m$-logarithm given in \cite[\S
12]{hain-macp}. Since Bloch-Wigner functions are single valued,
many of the technical problems one encounters when working with
multivalued are not present.

Our first result is:
\medskip

\noindent{\bf Theorem A.} {\it For each $m \le 4$ there exists
a real Grassmann $m$-logarithm.}
\medskip

The real Grassmann trilogarithm was first constructed as a
real Deligne cohomology class in \cite{yang:thesis}.
The symmetric group $\Sigma_{m+n+1}$ acts on $G^m_n$. If we insist
that each component of a representative of $L_m$ span a copy of
the alternating representation, then $L_m$ is unique.

A {\it generic real Grassmann $m$-logarithm} is an element of
$H_\bw^{2m}(U^m_\dot,\R(m))$, where $U^m_\dot$ is a Zariski open
subset of $G^m_\dot$ satisfying $U^m_0 = G^m_0$, whose image under
$$
\delta : H^{2m}_\bw(U^m\subdot,\R(m)) \to \Omega^m(G^m_0)
$$
is the volume form $\vol_m$. Observe that every real Grassmann
$m$-logarithm is a generic real Grassmann $m$-logarithm.

In the general case, we prove the following result, which is the
analogue for single valued polylogarithms of the main result of
\cite{hain:generic}.
\medskip

\noindent{\bf Theorem B.} {\it For each $m \ge 1$ there exists a
canonical generic real Grassmann $m$-logarithm.}
\medskip

The truncated simplicial variety $G^m_\dot$ can be viewed as a
``quotient'' of the classifying space $\bglm$ of rank $m$ vector
bundles.  One of the key points in the construction of real Grassmann
polylogarithms and in relating them to Chern classes is to show that
the ``alternating part'' of the universal Chern class
$$
c_m \in \Hd^{2m}(\bglm,\R(m))
$$
descends, along the ``quotient map'', to a class in
$\Hd^{2m}(G^m_\dot,\R(m))$. This approach is hinted at in the
survey \cite[p.~107]{bryl-zucker} of Brylinski and Zucker, although
the existence of this descended class is not evident. The descent
of the universal Chern class established in Section \ref{descent}.

The part of the cocycle of a (generic) real Grassmann $m$-logarithm
is a Bloch-Wigner function $D_m$ defined (generically) on $G^m_{m-1}$
which satisfies the $(2m+1)$-term functional equation
\begin{equation}\label{funct_eqn}
\sum_{j=0}^{2m}(-1)^j A_j^\ast D_m = 0.
\end{equation}
The function $D_1$ is simply $\log|\blank|$, the second function
$D_2$ is
the pullback of the Bloch-Wigner dilogarithm along the ``cross-ratio
map''
$G^2_1 \to \P^1 - \{0,1,\infty\}$ (cf. \cite[p.~403]{hain-macp}). The
functional equation (\ref{funct_eqn}) is the standard 5-term
equation.
The function $D_3$ is the single-valued trilogarithm whose existence
was established in \cite{hain-macp} and for which Goncharov
remarkably
expressed in terms of the single-valued classical trilogarithm in
\cite{goncharov:trilog}.

Recall that the {\it rank} filtration
$$
0 = r_0 K_m(R) \subseteq r_1 K_m(R) \subseteq r_2 K_m(R) \subseteq
\dots
\subseteq K_m(R)_\Q := K_m(R)\otimes \Q
$$
of the rational $K$-groups of a ring $R$ is defined by
$$
r_j K_m(R) = \im \{H_m(GL_j(R),\Q)\rightarrow
H_m(GL(R),\Q)\}\cap K_m(R)_\Q.
$$
Here $K_m(R)_\Q$ is identified with its image in $H_m(GL(R),\Q)$
under
the Hurewicz homomorphism.  By Suslin's stability theorem
\cite{suslin},
$r_mK_m(F) = K_m(F)_\Q$ whenever $F$ is an infinite field.

In Section \ref{homology} we show that a Bloch-Wigner function
defined
generically on $G^m_{m-1}$ and which satisfies the functional
equation
(\ref{funct_eqn}) defines an element of
$$
H^{2m-1}(GL_m(\C)^\delta,\R),
$$
where $GL_m(\C)^\delta$ denotes the general linear group with the
discrete
topology.
Such a function therefore defines a mapping
$$
r_mK_{2m-1}(\spec \C) \to \R.
$$
\medskip

\noindent{\bf Theorem C.} {\it If $D_m$ is a Bloch-Wigner function
defined generically on $G^m_{m-1}$ associated to the canonical choice
of a generic real Grassmann $m$-logarithm, then the associated map
$r_mK_{2m-1}(\C) \to \R$ is equals the restriction of the
Beilinson-Chern class
$$
c_m^B : K_{2m-1}(\spec\C) \to \Hd^1(\spec\C,\R(m))\cong \R
$$
to $r_mK_{2m-1}(\spec\C)$.}
\medskip

If $k$ is a number field then, by \cite{yang:rank,borel-yang},
$$
r_mK_{2m-1}(\spec k)=K_{2m-1}(\spec k)
$$
It follows that the {\it regulator mapping}
$$
c_p^B : K_{2m-1}(\spec k) \to \Hd^1(\spec k,\R(m)) \approx \R^{d_m},
$$
where $d_m$ is $r_1 + r_2$ or $r_2$ according to whether $m$ is odd
or even,
can be expressed in terms of  a generic real $m$-logarithm function.
It then follows by standard arguments, using fundamental results
of Borel, that when $m>1$, the value $\zeta_k(m)$ of the Dedekind
zeta
function of $k$ can be expressed, up to the product of a suitable
power
of $\pi$, a non-zero rational number, and the square root of the
discriminant of $k$, as a determinant of values of $D_m$ evaluated
at certain $k$-rational points of $G^m_{m-1}$. This generalizes the
classical theorem of Dedekind for the residue at $s=1$ of
$\zeta_k(s)$
and similar formulas for the values at $s=2$ due to Bloch and Suslin,
and at $s=3$ due to Goncharov \cite{goncharov:trilog} and Yang
\cite{yang:annoc,yang:thesis}.

An ultimate goal is to express a generic real Bloch-Wigner
$m$-logarithm
function in terms of the single-valued classical
$m$-logarithm. In this case, the Borel regulator would be expressed
in
terms of a determinant of values of the single-valued classical
$m$-logarithm at $k$ rational points of
$\P^1-\{0,1,\infty\}$.\footnote{
This is a weak statement of Zagier's conjecture.} More importantly,
it would show how to use the single-valued classical $m$-logarithm to
define the $m$th regulator. To date this has only
been done when $m\le3$: $m=2$ by Bloch and Suslin \cite{bloch}, $m=3$
by Goncharov \cite{goncharov:trilog}.

More generally, a generic real Grassmann $m$-logarithm defines a
function
$$
r_mK_n(\eta_X) \to \Hd^{2m-n}(\eta_X,\R(m))
$$
where $\eta_X$ denotes the generic point of the complex algebraic
variety
$X$.
\medskip

\noindent{\bf Theorem D.} {\it The function
$$
r_mK_n(\eta_X) \to \Hd^{2m-n}(\eta_X,\R(m))
$$
associated to the canonical choice of a generic real Grassmann
polylogarithm is the restriction of the
Beilinson-Chern class $c_m^B$ to $r_mK_n(\eta_X)$.}

\begin{remark} One should be able to write down the Chern
class on all of $K\subdot(\eta_X)$ using Goncharov's work
\cite{goncharov:chern} by
constructing ``generic bi-Grassmann polylogarithms'', but we have not
yet done this.
\end{remark}

\noindent{\it Conventions.} In this paper, all simplicial objects are
strict --- that is, they are functors from the category $\Delta$ of
finite ordinals and {\it strictly} order preserving maps to, say,
the category of algebraic varieties.

As is standard, the finite set $\{0,1,\dots,n\}$ with its natural
ordering will be denoted by $[n]$. Let $r$ and $s$ be positive
integers
with $r \le s$. Denote the full subcategory of $\Delta$ whose objects
are the ordinals $[n]$ with $r\le n \le s$ by $\Delta[r,s]$. An
$(r,s)$
{\it truncated} simplicial object of a category $\calC$ is a
contravariant
functor from $\Delta[r,s]$ to $\calC$.

The word {\it simplicial} will be used generically to refer to both
simplicial objects and truncated simplicial objects. However, we will
use the word truncated when we do want to emphasize the difference.
The
distinction will be made in \S \ref{descent} where is will be
significant.

By Deligne cohomology, we shall mean Beilinson's refined version
of Deligne cohomology as defined in \cite{beilinson:hodge} which is
sometimes called absolute Hodge cohomology. It can
be expressed as an extension
$$
0 \to \Ext^1_\H(\Q,H^{k-1}(X,\Q(p))) \to \Hd^k(X,\Q(p)) \to
\Hom_\H(\Q,H^k(X,\Q(p)))\to 0
$$
where $\H$ denotes the category of $\Q$ mixed Hodge structures.

To avoid confusion between, say, the $K$-theory of the ring $\C$ and
the variety $\C$, we shall view $K$ as a functor on schemes. We shall
therefore denote the $K$-theory of a ring $R$ by $K_\dot(\spec R)$.

\section{Bloch-Wigner Forms}

In this section, we introduce the complex of Bloch-Wigner forms on
a smooth complex algebraic variety.

Let $X$ be a smooth variety. Choose any smooth compactification
$\Xbar$ of $X$ where $D:=\Xbar - X$ is a normal crossings divisor.
Denote the complex of global meromorphic forms on $\Xbar$ which
are holomorphic on $X$ and have logarithmic singularities along $D$
by $\Omega^\bullet(\Xbar \log D)$.  This maps injectively to the
complex of holomorphic forms on $X$. Its image does not depend on
the choice of $\Xbar$, \cite[(3.3)]{hain-macp}, and will be denoted
by $\Omega^\bullet(X)$.

We now recall the definition of Bloch-Wigner functions on $X$ from
\cite[\S 11]{hain-macp}. Denote the algebra of all iterated
integrals of elements of $\Omega^1(X)$ by $A(X)$, and those that are
relatively closed by $H^0(A(X))$. Fix a base point $x\in X$. Taking a
relatively closed iterated integral
$$
\sum \int \omega_{i_1}\omega_{i_2}\dots\omega_{i_r}
$$
to the function
$$
z \mapsto \int_x^z \omega_{i_1}\omega_{i_2}\dots\omega_{i_r}
$$
defines an injective algebra homomorphism
$$
H^0(A(X))\rightarrow \widetilde{E}^0(X,x)
$$
where $\widetilde{E}^0(X,x)$ denotes the multivalued differentiable
functions on $X$ (see \cite[\S 2]{hain-macp}). The image of the
above map will be denoted by $\otilde(X,x)$.

Let $\otilde_\R(X,x)$ denote the subalgebra of the algebra of
multivalued, real valued functions on $X$ generated (as an
$\R$-algebra) by the real and imaginary parts of elements of
$\otilde(X,x)$. The algebra $\bw(X)$ is defined to be the subalgebra
of $\otilde_\R(X,x)$ consisting of single valued functions.
Equivalently, $\bw(X)$ is the subalgebra of
$\otilde_\R(X,x)$ invariant under monodromy:
$$
\bw(X)=\otilde_\R(X,x)^{\pi_1(X,x)}.
$$
Although this construction makes use of a base point, the ring
$\bw(X)$ itself depends only on $X$.  The assignment of $\bw(X)$
to $X$ is a contravariant functor from the category of smooth
complex algebraic varieties to the category of $\R$-algebras.
We call $\bw(X)$ {\it the ring of Bloch-Wigner functions on $X$}.

There is a natural weight filtration on $\Omega^\bullet(X)$. It
induces one on $A(X)$ by linear algebra,
and one on $H^0(A(X))$ by restriction.  This weight filtration passes
to $\otilde(X,x)$ and eventually to a weight filtration on
$\bw(X)$. This filtration is independent of all choices (cf.
\cite[\S 11]{hain-macp}).

Denote the subalgebra of the real de~Rham complex of $X$ generated by
the real and imaginary parts of elements of $\Omega^\dot(X)$
by $\Omega^\dot_\R(X)$.

\begin{definition}
The {\it complex of Bloch-Wigner forms on $X$\/} is defined to
be the sub-algebra of the de~Rham complex
$$
\obw^\bullet(X) = \bw(X)\cdot\Omega^\bullet_\R(X),
$$
generated by $\bw(X)$ and $\Omega^\bullet_\R(X)$.
\end{definition}

It is not difficult to see that the image of the natural inclusion
$$
\obw^\bullet (X)\hookrightarrow E^\bullet(X)
$$
is closed under $d$, so that $\obw^\dot(X)$ is a d.g.\ algebra.

The weight filtration of $\Omega^\bullet(X)$ induces a weight
filtration
on $\Omega_\R^\bullet(X)$. Taking the convolution of the weight
filtrations
of $\bw(X)$ and $\Omega^\bullet(X)$ we obtain a natural weight
filtration on $\obw^\bullet(X)$.

\section{Bloch-Wigner cohomology}
\label{bw-coho}

We introduce a natural analogue of Deligne (or more accurately,
absolute Hodge) cohomology with coefficients in $\R(m)$ which is
defined using Bloch-Wigner forms.

For a real vector space $V$, we denote $V\otimes \R(m)$ by $V(m)$,
and we identify the quotient $V_\C/V(m)$ with $V(m-1)$ using the
natural projection associated to the decomposition
$V_\C = V_\R(m-1) \oplus V_\R(m)$.

Suppose that $X$ is a smooth complex algebraic manifold.
Taking the value of an element of $\Omega^\bullet(X)$ mod $\R(m)$
defines a map
$$
\Omega^\bullet(X) \to \Omega^\bullet_\R(m-1)
$$
which preserves the weight filtrations. Composing with the
canonical inclusion $\Omega^\bullet_\R(X) \hookrightarrow
\obw\supdot(X)$ twisted by $\R(m-1)$, we obtain a natural weight
filtration preserving map
$$
\Omega\supdot(X) \rightarrow \obw\supdot(X)(m-1).
$$
For each $m\ge 0$, define a complex $\R(m)_\bw\supdot(X)$ by
$$
\R(m)_\bw\supdot(X) := \cone[F^mW_{2m}\Omega\supdot(X)
\rightarrow W_{2m}\obw\supdot(X)(m-1)][-1].
$$
Define $H\supdot_\bw(X,\R(m))$ to be the cohomology of this complex.
We shall call it the {\it $\bw$-cohomology of $X$ with coefficients
in $\R(m)$}. It is clearly functorial in $X$.

Denote by $\Hd^\bullet(X,\R(m))$ the absolute Hodge cohomology of
$X$ with coefficients in $\R(m)$.  This is Beilinson's refined
version
of Deligne-Beilinson cohomology \cite{beilinson:hodge}.
It is defined as the cohomology of the complex
$$
\R(m)_{\cal D}^\bullet(X) := \cone[F^mW_{2m}A_\C\supdot(X)
\rightarrow W_{2m}A_\R\supdot(X)(m-1)][-1],
$$
where
$$
A = ((A_\R^\bullet,W_\bullet),(A_\C^\bullet,W_\bullet, F^\bullet))
$$
is a real mixed Hodge complex for $X$. (Note that the weight
filtration used is the {\it filtration decal\'ee} of Deligne
\cite[p.~15]{deligne:II}.)

\begin{proposition}
There is a natural map $H^\dot_\bw(X,\R(m)) \to \Hd^\dot(X,\R(m))$.
\end{proposition}

\begin{pf}
We use the real mixed Hodge complex described in
\cite[p.~73]{durfee-hain}. Denote it by
$$
\left( (\Atilde_\R^\dot(X),W_\dot),
(\Atilde_\C^\dot(X),W_\dot,F^\dot)\right).
$$
The mixed Hodge complex with the {\it filtration decal\'ee} used to
compute Deligne cohomology is
$$
\left( (A_\R^\dot(X),W_\dot),(A_\C^\dot(X),W_\dot,F^\dot)\right),
$$
where
$$
W_lA^k = \left\{ a \in W_{l-k}\Atilde^k : da
\in W_{l-k-1}\Atilde^{k+1}\right\}
$$
and
$$
A^\dot = \bigcup_{l \ge 0} W_lA^\dot
$$
(cf. \cite[p.~145]{morgan}.)
It is straightforward to show that the image of the inclusion of
$\Omega_\bw^\dot(X)$ into $\Atilde_\R^\dot(X)$ is contained in
$A^\dot(X)$ and that the corresponding map
$$
\Omega_\bw^\dot(X) \to A_\R^\dot(X)
$$
preserves $W_\dot$. The result follows as the diagram
$$
\begin{matrix}
W_{2m}\Omega^\dot(X) & \to & W_{2m}\Omega^\dot_\bw(X)(m-1) \cr
\downarrow & & \downarrow \cr
W_{2m}A^\dot_\C(X) & \to & W_{2m}A^\dot_\R(X)(m-1)\cr
\end{matrix}
$$
in which the horizontal maps are reduction of values mod $\R(m)$,
commutes.
\end{pf}

We shall need the simplicial analogue of this result. The proof is
similar to that of the previous result.

\begin{proposition}
For each smooth simplicial variety $X_\dot$, there is a natural
map
$$
H_\bw^\dot(X_\dot,\R(m)) \to \Hd^\dot(X_\dot,\R(m)). \qed
$$
\end{proposition}

Our polylogarithms lie in the $\bw$-cohomology of certain
simplicial varieties. In order to construct polylogarithms, we
shall need to compare the $\bw$- and Deligne cohomologies of
these simplicial varieties.

Recall that the irregularity $q(X)$ of a smooth complex algebraic
variety $X$ is defined
by
$$
q(X) := \dim W_1H^1(X;\Q)/2 = \dim H^{1,0}(\Xbar),
$$
where $\Xbar$ is any smooth completion of $X$. Most varieties in this
paper will satisfy $q(X)=0$. This condition is satisfied by all
Zariski
open subsets of simply connected smooth varieties. In particular, it
is
satisfied by each $G^m_n(\C)$.

Suppose that $X$ is a topological space with finitely  generated
fundamental group. Denote the $\C$-form of the Malcev Lie algebra
associated to $\pi_1(X,x)$ by $\g(X,x)$. This is a topological Lie
algebra. There is a canonical homomorphism
$$
\Hcts^\dot(\g(X,x)) \to H^\dot(X,\C),
$$
induced, for example, by the homomorphism
$$
\Lambda^\dot(\g^\ast)\rightarrow \Omega^\dot(X)\hookrightarrow
E^\dot(X),
$$
(see \cite[\S 7]{hain-macp}). Recall from \cite[(8.3)]{hain-macp}
that
$X$ is a rational $n$-$K(\pi,1)$ if this map is an isomorphism in
degrees $\le n$ and injective in degree $n+1$.

\begin{theorem}\label{summand}
Suppose that $X_\dot$ is a simplicial variety where each $X_m$ is
smooth and has $q=0$. If, for each $m$, $X_m$ is a rational $(n-m)$-%
$K(\pi,1)$, then the natural map
$$
H_\bw^t(X_\dot,\R(m)) \to \Hd^t(X_\dot,\R(m))
$$
has a canonical splitting whenever $t\le n$. In particular, this
map is surjective in degrees $\le n$.
\end{theorem}

\section{Simplicial Spaces with Symmetric Group Actions}
\label{S-variety}

Certain simplicial spaces come equipped with actions of symmetric
groups
on their spaces of simplices. This symmetric group action gives an
extra
algebraic structure to the cohomology of complexes associated to such
simplicial spaces. This was first explored in \cite[\S9]{hain-macp}.
Here
we develop those ideas a little further.

\begin{definition}
A simplicial topological space $X\subdot$ is called a
$\Sigma\subdot$-space if there exits, on each $X_n$, a continuous
action of the symmetric group $\Sigma_{n+1}$ on $n+1$ letters which
is
compatible with the simplicial structure of $X\subdot$ in the
sense that the face maps
$$
A_j:X_{n+1}\rightarrow X_n,
$$
satisfy the conditions
\begin{equation*}
A_j\circ (i-1,i) = \begin{cases}
	             (i-1, i)\circ A_j & j> i;\\
		     A_{i-1} & j=i;\\
		     A_i  & j=i-1;\\
		     (i-2,i-1)\circ A_j & j< i.
		     \end{cases}
\end{equation*}
\end{definition}

When the simplicial space $X\subdot$ has extra
structure, (examples being when $X\subdot$ is a simplicial manifold
or
a simplicial variety), we require the $\Sigma\subdot$ action to
preserve
this additional structure.

\begin{example}
Let $X$ be an arbitrary topological
space. Consider the simplicial space $X_\bullet$ where
$X_n=X^{n+1}$ with the obvious face
maps. The symmetric group $\Sigma_{n+1}$ acts on $X_n$ by
permuting the factors. It is routine to check that $X\subdot$ is a
$\Sigma_\bullet$-space.
\end{example}

\begin{example}
This type of simplicial space arises in the construction of
classifying spaces of principal $G$ bundles.
Let $X$ be a principal $G$-space, where $G$ is a topological group.
Define a simplicial space $(X/G)_\bullet$ by letting
$$
(X/G)_n= X^{n+1}/G,
$$
where $G$ acts diagonally of $X^{n+1}$. Since the
$\Sigma_{n+1}$-action,
permuting the factors of $X^{n+1}$, commutes with the $G$-action,
$(X/G)_n$ inherits
a $\Sigma_{n+1}$-action which gives it the structure of a
$\Sigma_\bullet$-space. Note that
$(X/G)_\bullet$ is just the standard simplicial model of the
classifying
space $B_\bullet G$ of principal $G$-bundles. When $X$ is an
algebraic
variety and $G$ is an algebraic group which acts on $X$ algebraically
and where each $X^n/G$ is an algebraic variety, $(X/G)_\bullet$ is a
simplicial variety.
\end{example}

We now describe several $\Sigma_\dot$-varieties to be used in the
paper.

The scheme $G^m_n$ introduced in the introduction has an alternative
description (see \cite[(5.6)]{hain-macp}). We say that $m+n$ vectors
in
a vector space $k^m$ are in {\it general position} if each $m$ of
them are
linearly independent. The alternative description of $G^m_n$ is:
$$
G^m_n(k)=\{(n+m+1)\text{-tuples}\; (v_0,\dots,v_{m+n}) \;
\text{in $k^m$ in general position}\}/GL_m(k).
$$
We will denote the point of $G^m_n$ corresponding to the orbit of
$(v_0,\dots, v_{m+n})$ by $[v_0,\dots, v_{m+n}]$.
The face map
$$
A_i: G^{m+1}_n  \rightarrow  G^m_n
$$
is defined by
$$
[v_0,\dots,v_{m+n+1}] \mapsto [v_0,\dots,\hat{v}_i,\dots,v_{m+n+1}].
$$
It is clear that $G^m_n$ has a $\Sigma_{m+n+1}$ action. If we place
$G^m_n$ in dimension $n+m$, then the $G^m_n$ with $0\le n \le m$,
together with the face maps form an $(m,2m)$ truncated
$\Sigma\subdot$-variety $G^m\subdot$.

It is convenient to complete the truncated variety $G^m\subdot$ to a
simplicial variety $\bgm$ by adding a point in each degree less
than $m$. Observe that $\bgm$ is also a $\Sigma\subdot$-variety.
Similarly, each Zariski open subset $U^m_\dot$ of $G^m_\dot$ can be
completed to a simplicial variety $\bum$. If $U^m\subdot$ is
invariant
under the symmetric group actions on $G^m\subdot$, then $\bum$ is
also a $\Sigma\subdot$-variety.

Define
$$
E_nG^m(k) = \{(n+m+1)\text{-tuples}\; (v_0,\dots,v_{m+n}) \;
\text{in $k^m$ in general position}\},
$$
and let $\egm$ be the set of $E_nG^m$ with $n\ge 0$. Then with the
obvious
face maps, with the obvious face maps, $\egm$ forms a
$\Sigma_\dot$-variety.
The group scheme $GL_m$ acts on $\egm$ via the diagonal action.
Observe that the quotient is $\bgm$. Denote the projection map by
$$
\pi_m: \egm\rightarrow \bgm.
$$

For a Zariski open subset $U^m_\dot$ of $G^m_\dot$, let
$\eum = \pi_m^{-1}(\bum)$. Then $\eum$ is a simplicial open
subvariety of $\egm$, which is a $\Sigma\subdot$-variety
when $U^m_\dot$ is.

Let $V$ be a $k$-module with a $\Sigma_n$-action, where $k$ is a
field
of characteristic 0. Define the {\it alternating operator}
$$
\Alt_n: V\rightarrow V
$$
by
\begin{equation*}
\Alt_n(v) = \frac{1}{n!}\sum_{\sigma\in
\Sigma_n}\sgn(\sigma)\sigma(v).
\end{equation*}

When $n>1$, an element $v\in V$ is called an {\it alternating element
of $V$} if $\Alt_n(v)=v$. Let $sV$ denote the submodule of $V$
consisting of all of its alternating elements. (This will be called
the {\it alternating part\/} of $V$ in the sequel.)  The {\it sign
decomposition\/} of $V$ is the decomposition
$$
V= sV\oplus rV
$$
where $rV= \ker\; \Alt_n$ is the unique $\Sigma_n$-invariant
complement of $sV$.  By convention, we define $sV = V$ and
$rV=0$ when $n=1$. This is necessary in order that the following
result
hold.

\begin{definition} (cf. \cite[(9.4)]{hain-macp})
A cosimplicial $k$-module $M\supdot$ ($k$ a field
of characteristic 0) is called a {\it $\Sigma_\bullet$-module\/} if
each
$M^n$ has a $\Sigma_{n+1}$-action, and the face maps
$$A_j: M^{n-1}\rightarrow M^n$$
satisfy the following relations
\begin{equation*}
(i-1,i)\circ A_j =
\begin{cases}
A_j\circ(i-2,i-1), & j<i-1;\\
A_i, & j= i-1;\\
A_{i-1}, & j=i;\\
A_j\circ(i-1,i) & j> i,
\end{cases}
\end{equation*}
for $ j=0,\dots,m$ and $i=1,\dots, m-1$.
\end{definition}

Natural examples of cosimplicial $\Sigma_\bullet$-modules can be
obtained
by applying a contravariant $k$-module valued functor to a
$\Sigma_\bullet$-space.

As usual, let us define $A^\ast: M^{n-1}\rightarrow M^n$ to be
$\sum^n_{j=0}(-1)^j A_j$. The sign decomposition generalizes to
$\Sigma_\bullet$-cosimplicial modules. The following result is proved
in \cite[(9.5)]{hain-macp}.

\begin{lemma}\label{signd}
The differential $A^\ast$ preserves the sign decomposition
$$
M^l = sM^l\oplus rM^l,
$$
for $l=0,1,2,\dots$. In particular, $A^\ast sM^n$ is a
$\Sigma_{n+2}$-submodule of $M^{n+1}$. \qed
\end{lemma}

Define the cohomology of a cosimplicial module $M^\bullet$ by
$$
H^n(M^\bullet)=\ker(A^\ast: M^n\rightarrow
M^{n+1})/\im(A^\ast: M^{n-1}\rightarrow M^n).
$$
The following corollary is an immediate consequence of the previous
lemma.

\begin{corollary}
The cohomology groups of a cosimplicial $\Sigma_\bullet$-module
$M^\bullet$ have a sign decomposition
$$
H^\bullet(M^\bullet)=sH^\bullet(M^\bullet)\oplus
rH^\bullet(M^\bullet).
$$
The decomposition is natural with respect to
$\Sigma_\bullet$-invariant maps between cosimplicial
$\Sigma_\bullet$-modules.
\end{corollary}

Applying the de~Rham complex functor, the Deligne-Beilinson cochain
complex functor, or the $\bw$-cochain complex functor to a
$\Sigma_\bullet$-variety, we obtain natural
examples of $\Sigma_\bullet$-cosimplicial modules. The following
result
is an immediate consequence of the previous result.

\begin{theorem}
If $X_\bullet$ is a smooth complex $\Sigma_\bullet$-variety, then
the de~Rham cohomology, Deligne-Beilinson cohomology and the $\bw$-%
cohomology of $X_\bullet$ have sign decompositions
$$
H^\bullet(X_\bullet)=sH^\bullet(X_\bullet)\oplus
rH^\bullet(X_\bullet),
$$
$$
\Hd^\bullet(X_\bullet,\Lambda(m))=
s\Hd^\bullet(X_\bullet,\Lambda(m))\oplus
r\Hd^\bullet(X_\bullet,\Lambda(m)),
$$
$$
H_\bw^\bullet(X_\bullet,\R(m))= sH_\bw^\bullet(X_\bullet,\R(m))\oplus
rH_\bw^\bullet(X_\bullet,\R(m)),
$$
which are all natural with respect to $\Sigma_\bullet$-invariant map
between smooth $\Sigma_\bullet$-varieties. Moreover, the natural
maps
$$
H^\dot_\bw(X_\dot,\R(m)) \to \Hd^\dot(X_\dot,\R(m)) \to
H^\dot(X_\dot,\R(m))
$$
each preserve the sign decomposition. \qed
\end{theorem}

\section{Real Grassmann Polylogarithms}
\label{polylog}

In this section we shall view the truncated simplicial variety
$G^m_\dot$ as a $\Sigma_\dot$-variety with the $\Sigma_\dot$
action described in \S \ref{S-variety}.

Denote the coordinates of $\P^m$ by $[x_0,x_1,\dots,x_m]$. Denote
the hyperplane $x_j = 0$ by $H_j$. There is a natural identification
of $G^m_0$ with $\P^m - \cup_{j=0}^m H_j$. This can be identified
with
$(\C^\ast)^m$ by identifying $(x_1,\dots,x_m)\in (\C^\ast)^m$ with
$[1,x_1,\dots,x_m]\in \P^m$. Set
$$
\vol_m = \frac{dx_1}{x_1}\wedge\dots \wedge \frac{dx_m}{x_m}.
$$
This is an element of
$$
s\Omega^m(G^m_0) \cong sH^m(G^m_0,\C).
$$

Recall that $G^m_n$ is placed in dimension $m+n$. There is a
canonical
homomorphism
$$
s\Hd^{2m}(G^m_\dot,\R(m)) \to sH^m(G^m_0,\R(m-1))
$$
induced by the inclusion $G^m_0 \hookrightarrow G^m\subdot$.

\begin{definition}
A {\it real Grassmann $m$-logarithm} is an element $L_m$ of
$$
sH_\bw^{2m}(G^m_\dot,\R(m))
$$
whose image under the composite
$$
sH_\bw^{2m}(G^m_\dot,\R(m)) \to s\Hd^{2m}(G^m_\dot,\R(m)) \to
sH^m(G^m_0,\C)
$$
is $\vol_m$. A cocycle in $s\R_\bw^\dot(G^m_\dot)(m)$ that
represents $L_m$ will be called a {\it real Grassmann $m$-cocycle.}
Finally, the part of a real Grassmann $m$-cocycle that lies in
$sW_{2m}\bw(G^m_{m-1})$ will be called a {\it real Grassmann
$m$-logarithm function}.
\end{definition}

Observe that a real Grassmann $m$-logarithm $D_m$ satisfies the
$(2m+1)$-term functional equation
$$
A^\ast D_m := \sum_{j=0}^{2m} (-1)^j A_j^\ast D_m = 0,
$$
where $A_j : G^m_m \to G^m_{m-1}$, $j=0,\dots,2m$, are the face maps,
as well as the skew symmetry property
$$
\sigma^\ast D_m = \sgn(\sigma) D_m
$$
for all $\sigma \in \Sigma_{2m}$.

Now suppose that $U^m_\dot$ is a Zariski open subset of $G^m_\dot$.
(That is, for each $n$, $U^m_n$ is a Zariski open subset of $G^m_n$
and the inclusion $U^m_\dot \hookrightarrow G^m_\dot$ is a simplicial
map.) Suppose further that $U^m_0 = G^m_0$, that $U^m_\dot$ is
mapped into itself under the action of the symmetric groups on
$G^m_\dot$, and that the condition
\begin{equation}\label{condition}
\text{\it Each fiber of each face map $A_j : U^m_l \to U^m_{l-1}$ is
non-empty}
\end{equation}
is satisfied.

As above, the inclusion $G^m_0 \hookrightarrow U^m_\dot$
induces a canonical homomorphism
$$
s\Hd^{2m}(U^m_\dot,\R(m)) \to sH^m(G^m_0,\R(m-1)).
$$

\begin{definition}\label{hain:generic}
A {\it generic real Grassmann $m$-logarithm} is an element $L_m$ of
$$
sH_\bw^{2m}(U^m_\dot,\R(m)),
$$
where $U^m_\dot$ is a Zariski open
subvariety of $G^p_\dot$ that satisfies the conditions in the
previous paragraph, whose image under the composite
$$
sH_\bw^{2m}(U^m_\dot,\R(m)) \to s\Hd^{2m}(U^m_\dot,\R(m)) \to
sH^m(G^m_0,\C)
$$
is $\vol_m$. A cocycle in $s\R_\bw^\dot(U^p_\dot)(m)$ that
represents $L_p$ will be called a {\it generic real Grassmann
$m$-cocycle.}
Finally, the part of a generic real Grassmann $m$-cocycle that lies
in
$sW_{2m}\bw(U^m_{m-1})$ will be called a {\it generic real Grassmann
$m$
logarithm function}.
\end{definition}

Observe that every real Grassmann $m$-logarithm (resp.\ cocycle,
function)
is a generic real Grassmann $m$-logarithm (resp.\ cocycle, function).

As in the case of a Grassmann $m$-logarithm, a generic Grassmann
$m$-logarithm $D_m$ satisfies the $(2m+1)$-term functional equation
$$
A^\ast D_m = 0
$$
and the symmetry relation
$$
\sigma^\ast D_m = \sgn (\sigma) D_m
$$
for each $\sigma \in \Sigma_{2m}$.

In the next section, we will prove that a generic real Grassmann
$m$-logarithm (and therefore, every real Grassmann $m$-logarithm)
defines a cohomology class
$$
d_m \in H^{2m-1}(GL_m(\C),\R(m-1)),
$$
and in the succeeding section, that $D_m$ defines a cohomology class
in the continuous cohomology
$$
\delta_m \in \Hcts^{2m-1}(GL_m(\C),\R(m-1))
$$
such that the image of $\delta_m$ in $H^{2m-1}(GL_m(\C),\R(m-1))$
is $d_m$.

\section{Up to the 3-log}\label{up}

In this section, we establish the existence of the
first 3 real Grassmann logarithms.  This results is more or less
known from \cite{bloch}, \cite{hain-macp}, \cite{yang:thesis} and
\cite{goncharov:trilog}.

\begin{proposition}
If $m\le 3$, then $sW_{2m}H^{2m}(\bgm,\Q)$ has dimension one
and is spanned by the volume form $\vol_m$, while
$sW_{2m-1}H^{2m-1}(\bgm,\Q)$ is trivial.
\end{proposition}

\begin{pf}
When $m\le 3$, the stronger result for all the cohomology (rather
than just the alternating part) was proved by direct computation
by Hain and MacPherson (cf. \cite[(12.6)]{hain-macp}), although the
details were not given. Here we  give a complete proof of the weaker
statement given in the proposition.

The point is that when $m\le 3$, $G^m_n$ is a rational
$(m-n)$-$K(\pi,1)$
for all $n$. This is proved in \cite[\S 8]{hain-macp}. This implies
that for such $m$ and $n$, the cup product
$$
\Lambda^k H^1(G^m_n,\Q) \to H^k(G^m_n,\Q)
$$
is surjective, provided that $k\le m-n$. One can easily show that in
these cases, $s\Lambda^k H^1(G^m_n,\Q)=0$, except when $k=m$ and
$n=0$,
in which case
$$
s\Lambda^m H^1(G^m_0,\Q) \cong sH^m(G^m_0,\Q).
$$
The result now follows from the fact that the standard spectral
sequence
that converges to $H^\dot(G^m_\dot,\Q)$ is compatible with the
$r\oplus s$
decomposition.
\end{pf}

As a corollary, we obtain the existence and uniqueness of real
Grassmann $m$-logarithms for $m=1,2,3,4$.

\begin{corollary}\label{one-dim}
If $m\le 3$, then the natural map
$$
s\Hd^{2m}(\bgm,\R(m)) \to \Hd^{2m}(G^m_0,\R(m))\cong \R\vol_m
$$
is an isomorphism.
\end{corollary}

\begin{pf}
The proposition follows immediately from the following short exact
sequence
\begin{multline*}
0 \to \Ext^1_\H(\Q,H^{2m-1}(X,\Q(m))) \to \\
\Hd^{2m}(X,\Q(m)) \to \Hom_\H(\Q,H^{2m}(X,\Q(m)))\to 0
\end{multline*}
since by the previous proposition
$$
\Ext^1_\H(\Q,H^{2m-1}(X,\Q(m)))=0
$$
and
$$
\Hom_\H(\Q,H^{2m}(X,\Q(m)))
$$
is one-dimensional and generated by the volume form $\vol_m$.
\end{pf}

\begin{corollary}
If $m \le 3$, there is a canonical $m$-logarithm.
\end{corollary}

\begin{pf}
Since in each case $G^m_n$ is a rational $(m-n)$-$K(\pi,1)$, there is
a canonical splitting of the map
$$
sH^{2m}_\bw(G^m_\dot,\R(m)) \to s\Hd^{2m}(G^m_\dot,\R(m)) \cong
\R\vol_m
$$
by (\ref{summand}). The $m$-logarithm is the image of $\vol_m$ under
this
splitting.
\end{pf}

\section{Homology of $GL_m$}
\label{homology}

Let $k$ be an infinite field. In this section we show that a function
$$
f: G^m_n(k) \rightarrow \R
$$
satisfying the functional equation
$$
A^\ast f = 0
$$
determines an element of $H^{m+n}(GL_m(k);\R)$.
In fact, with a little bit more work, we will show that the same
holds
even if $f$ is only defined on a Zariski open subvariety $U^m_n$ of
$G^m_n$.

Recall that for an abstract group $G$, the group cohomology of $G$
with
coefficients in an abelian group $V$ (viewed as a trivial $G$-module)
can be computed by taking the homology of the $G$ invariants of the
complex $C\supdot(G,V)$, where $C^n(G,V)$ is the group of functions
$$
f: \underbrace{G\times\dots\times G}_{n+1}\rightarrow V
$$
and the coboundary map is defined by
$$
\delta f(g_0,\dots,g_{n+1}) =
\sum_{i=0}^{n+1}(-1)^if(g_0,\dots,\hat{g}_i,\dots,g_{n+1}).
$$
The group $G$ acts on $C\supdot(G,A)$ on the right via the formula
$$
(f\cdot g)(g_0,\dots,g_n) := f(gg_0,\dots,gg_n).
$$
The $G$-invariants of the complex $C\supdot(G,A)$ will be denoted by
$C\supdot(G,V)^G$.

\begin{variant}\label{coho_def}
Recall that if $G$ is a topological group,
then the continuous group cohomology of $G$ is defined using
continuous functions $f: G^{n+1} \to V$ in place of arbitrary
functions
in the definition above. When $G$ is a Lie group, the locally-$L_p$
cohomology of $G$ is defined using cochains $f:G^{n+1} \to V$ that
are
locally-$L_p$ with respect to the measure on $G$ given by a left
invariant
volume form. The cochain complex of continuous cochains of $G$ will
be
denoted by $C_{\text{cts}}^\dot(G,V)$, and the locally-$L_p$ cochains
by
$C_{\text{loc-$L_p$}}^\dot(G,V)$. The continuous and locally-$L_p$
cohomology groups of a Lie group $G$ will be denoted
$H_{\text{cts}}^\dot(G,V)$ and $H_{\text{loc-$L_p$}}^\dot(G,V)$,
respectively.
\end{variant}
\medskip

Let $k$ be an extension field of $\Q$. Fix a non-zero vector $e$ in
$k^m$. Define the subset
$X^n_{G(k),e}$ of $E_nGL_m(k)$ to be\footnote{This is the discrete
analogue of the variety $B_nGL_m(\C)^\gen$ defined in \S
\ref{descent}.}
$$
\{(g_0,\dots,g_n)\in GL_m(k)^{n+1}
| \; (g_0e,\dots,g_ne)\in E_nG^m(k)\}.
$$
Denote the group of functions $X^n_{G(k),e}\rightarrow V$ by
$C^n_{G,e}(GL_m(k),V)$. Endowed with the boundary maps induced
from those of $C\supdot(GL_m(k),V)$, it is a complex. Denote the
subcomplex of $GL_m(k)$-invariants by
$C_{G,e}\supdot(GL_m(k),V)^{GL_m(k)}$.

\begin{proposition}
\label{resoln}
The chain map
$$
C\supdot(GL_m(k),V)^{GL_m(k)} \rightarrow
C\supdot_{G,e}(GL_m(k),V)^{GL_m(k)}
$$
induced via restriction from the natural chain map
$$ C\supdot(GL_m(k),V)\rightarrow C\supdot_{G,e}(GL_m(k),V) $$
induces an isomorphism on cohomology; i.e., there is a
natural isomorphism
$$ H^\bullet(C\supdot_{G,e}(GL_m(k),V)^{GL_m(k)})\cong
H^\bullet(GL_m(k),V). $$
\end{proposition}

\begin{pf} It suffices to prove that
$$
0\rightarrow V \rightarrow C\supdot_{G,e}(GL_m(k),V)
$$
is a resolution of $V$ by injective $GL_m(k)$-modules . To prove
this, we prove the dual statement.

Denote the tensor product of the free abelian group generated by the
points of  $X^m_{G(k),e}$ with $V$ by $C^{G,e}_n(GL_m(k),V)$. The
dual
statement is that
$$
0 \leftarrow V\leftarrow C\subdot^{G,e}(GL_m(k),V)
$$
is a projective resolution of $V$. Since each $C_n^{G,e}(GL_m(k),V)$
is a
free $GL_m(k)$-module, we need only establish exactness. Supposes
that
\begin{equation*}
\partial (\sum_l
a_k(g_{l,0},\dots,g_{l,n}))=\sum_k\sum_{i=0}^n(g_{l,0},\dots,
\hat{g}_{l,i},\dots,g_{l,n})=0,
\end{equation*}
where $(g_{l,0}e,\dots, g_{l,n}e)\in E_nG^m(k)$. By elementary linear
algebra, there exists $v\in k^m-\{0\}$ such that $(v,
g_{l,0}e,\dots, g_{l,n}e)\in E_{n+1}G^m(k)$ for each $l$. Pick
$g\in GL_m(k)$ such that $g e = v$. Then
$(g,g_{l,0},\dots,g_{l,n})$ lies in $X^{n+1}_{G,e}(k)$ for each $l$.
Now it is straightforward to check that
\begin{equation}
\label{exactness}
\partial(\sum_k
a_k(g,g_{k,0},\dots,g_{k,n}))=\sum_ka_k(g_{k,0},\dots,g_{k,n}).
\end{equation}
The exactness follows.
\end{pf}

Now suppose we are given a function
$$
f: G^m_n(k)\rightarrow V
$$
that satisfies $A^\ast f = 0$. Since $G^m_n(k)=B_{m+n}G^m(k)$ when
$n>0$, $f$ will induce a $GL_m(k)$-invariant map
$$
\tilde{f} : E_{m+n}G^m(k) \to V.
$$
As before, we fix a non-zero vector $e\in k^m$. Define a map
$$
f^e: X^{m+n}_{G,e}(k) \rightarrow V
$$
by
$$
f^e(g_0,\dots,g_{m+n})= \tilde{f}([g_0e,\dots,g_{m+n}e]).
$$
It is obvious that $f^e$ is a $GL_m(k)$-invariant
cocycle in $C^{m+n}_{G,e}(GL_m(k),V)$. In fact we have the
following result.

\begin{proposition}\label{groupclass}
The function $f^e$ represents a cohomology class in
$$
H^{m+n}(GL_m(k),V).
$$
Moreover, the cohomology class it represents is
independent of the choice of the base vector $e$ in $k^m-\{0\}$.
\end{proposition}

\begin{pf} The first statement follows immediately from
(\ref{resoln}).
To prove the second, suppose that $e'$ is another non-zero vector in
$k^m$. There exists a matrix $h\in GL_m(k)$ such that $e'= he$.
Define a chain map
$$
\phi_h: C\supdot(GL_m(k), V)^{GL_m(k)} \longrightarrow
C\supdot(GL_m(k), V)^{GL_m(k)}
$$
by
$$
\phi_h(f)(g_0,\dots, g_{m+n})= f(g_0h,\dots, g_{m+n}h),
$$
where $f\in C\supdot(GL_m(k),V)^GL_m(k)$. It is known (see, e.g.,
\cite[Chap. IV, Prop.~5.6]{maclane}) that this chain map induces the
identity map on
cohomology.

The map $\phi_h$ induces a chain map
$$
\phi^{e,e'}_h: C_{G,e}\supdot(GL_m(k), V)^{GL_m(k)}
\longrightarrow C_{G,e'}\supdot(GL_m(k), V)^{GL_m(k)},
$$
which carries $f^e$ to $f^{e'}$.
The proposition now follows from (\ref{resoln}) and the commutativity
of the following diagram.
$$
\begin{CD}
C\supdot(GL_m(k),V)^{GL_m(k)} @>\phi_h >>
C\supdot(GL_m(k),V)^{GL_m(k)}\\
@VVV                        @VVV \\
C\supdot_{G,e}(GL_m(k),V)^{GL_m(k)} @>\phi^{e,e'}_h >>
C\supdot_{G,e'}(GL_m(k),V)^{GL_m(k)}
\end{CD}
$$
\end{pf}

\begin{corollary}
\label{cor-dm}
If $D_m$ is a real Grassmann $m$-logarithm function, then $D_m$
defines an element of $H^{2m-1}(GL_m(\C),\R(m))$. \qed
\end{corollary}

\begin{remark}\label{grass_homo}
These results can be interpreted in terms of MacPherson's Grassmann
homology \cite{b-mcp-s}, as we shall now explain. Further discussion
of Grassmann homology and its relation to Suslin's work can be found
in \cite{wolf}.

Let $k$ be a field. Denote the free abelian group generated by the
points of $E_nG^m(k)$ by $C_n(k^m)$.
(This is denoted $C_n(GP(k^m))$ in \cite{suslin}.)
The face maps of $E_\dot G^m(k)$ induce a differential
$C_n(k^m) \to C_{n-1}(k^m)$. This is a resolution of the trivial
module \cite[Lemma 2.2]{suslin}.

The group $GL_m(k)$ acts diagonally on $C_n(k^m)$. Applying the
functor $\underline{\phantom{X}}\otimes_{GL_m(k)} k$, we
obtain the complex
$$
0\leftarrow C_0(k^m)_{GL_m(k)}\stackrel{d}{\leftarrow}
C_1(k^m)_{GL_m(k)}\stackrel{d}{\leftarrow} C_2(k^m)_{GL_m(k)}
\stackrel{d}{\leftarrow}\cdots
$$
where $C_n(k^m)_{GL_m(k)}$ denotes the $GL_m(k)$ coinvariants
$C_n(k^m)_{GL_m(k)}\otimes_G k$ --- $C_n(k^m)_{GL_m(k)}$ is simply
the free abelian group generated by the points of $B_nG^m(k)$.

The homology of this complex is called the {\it (extended) Grassmann
homology of $k$}, and is denoted by $GH_\bullet^m(\spec k)$. This
differs by a factor of $\Z$ from MacPherson's original definition in
dimension $m$ when $m$ is even.

The result (\ref{resoln}) simply says that there is a natural map
$$
H_n(GL_m(k)) \to GH^m_{m+n}(\spec k)
$$
and gives a formula for it. If $f : G^m_n(k) \to V$ satisfies the
functional equation $A^\ast f = 0$, then $f$ induces a map
$GH^m_{m+n}(\spec k) \to V$. The class $f^e$ of (\ref{groupclass}) is
simply the composite of these two maps.
\end{remark}

\begin{variant}\label{variant}
The corollary \ref{cor-dm} can be generalized to generic real
Grassmann
$m$-logarithm functions. Suppose that $U^m_\dot(k)$ is a Zariski open
subset of $G^m_n(k)$ that satisfies (\ref{condition}) in \S
\ref{polylog}. Now suppose that
that $f : U^m_n(k) \to V$ satisfies the
functional equation $A^\ast f = 0$. We will indicate briefly how to
modify the proof of (\ref{cor-dm}) to show that $f$ determines an
element $f^e$ of $H_{m+n}(GL_m(k),V)$.

Define
$$
E_{n+m}U^m = \{(v_0,\dots,v_n)|v_j\in k^m \text{ and }
[v_0,\dots,v_n] \in U^m_n(k)\}.
$$
Now replace $X^n_{G(k),e}$ by
$$
X^n_{U(k),e}:=\{(g_0,\dots,g_n)\in GL_m(k)^{n+1}
| \; (g_0e,\dots,g_ne)\in E_nU^m(k)\}.
$$
Let $\pi_0: E_{n+1}U^m(k)\rightarrow k^m-\{0\}$ denote the projection
map to the first component. Since $k$ is infinite, the condition
(\ref{condition}) of \S \ref{polylog} implies that for each finite
set
of points
$x_1, \dots, x_N$ of $E_nU^m(k)$,
$$
\bigcap_{l=1}^N \pi_0(A_0^{-1}(x_l))\subset k^m-\{0\}
$$
is non-empty. This condition is needed to prove that the complex
corresponding to the simplicial set $X^\dot_{U(k),e}$ is a resolution
of $\Z$.

Putting all this together, we have:
\end{variant}

\begin{corollary}
If $D_m$ is a generic real Grassmann $m$-logarithm function, then
$D_m$
defines an element of $H^{2m-1}(GL_m(\C),\R(m))$. \qed
\end{corollary}

We conclude this section with a useful
technical fact. Let $G$ be a discrete group.
The standard resolution for computing group cohomology comes from
the $\Sigma_\dot$-variety $BG$. It follows that the group homology
$H^\dot(G,\Q)$ has a sign decomposition.

\begin{proposition}\label{sign-copy}
The homology $H_\dot(G,\Q)$ consists entirely
of the alternating part. That is,
$rH^\bullet( C^\dot_\cts(B_\bullet G, V))= 0$.
\end{proposition}

\begin{pf}
Let $C_\dot(G,\Q) \stackrel{\epsilon}{\longrightarrow} \Q \to 0$
be the free resolution of the trivial module that comes from the
standard simplicial model of $BG$. Since $BG$ is a $\Sigma_\dot$-%
variety, this resolution has an $r\oplus s$ decomposition. Since
$sC_0(G,\Q) = C_0(G,\Q)$, and since $r$ and $s$ are exact functors,
it follows that
$$
sC_\dot(G,\Q) \stackrel{\epsilon}{\longrightarrow} \Q \to 0
$$
is a projective resolution of the trivial module. The result follows.
\end{pf}

The analogous results hold for both the continuous cohomology and
locally $L_p$  cohomology of a Lie  group. The proofs are similar.

\section{Bloch-Wigner Functions and Locally $L_p$ Cohomology}
\label{loc-l_p}

In this section, we show that a Bloch-Wigner function
$f:U^m_n(\C)\rightarrow \R$ that satisfies the
functional equation $A^\ast f=0$ represents a {\it continuous}
group cohomology class of $GL_m(\C)$ whose image in
$H_{m+n}(GL_m(\C),\R)$ is the class constructed from $f$ in Section
\ref{homology}.

Suppose that $M$ is an orientable smooth manifold of dimension
$n$ and that $\omega$ is a nowhere vanishing $n$-form on $M$. An
almost everywhere defined function $f$ on $M$ is said to be
{\it locally $L_p$} if each point $x$ of $M$ has a neighbourhood $U$
such that
$$
\int_U |f|^p \omega < \infty.
$$
This definition is independent of the choice of volume form
$\omega$.

\begin{lemma}\label{loc_pres}
Suppose that $\pi : M \to N$ is an orientation preserving proper map
between orientable manifolds of the same dimension. If
$f$ is an almost everywhere defined function on $N$ whose pullback is
almost everywhere defined on $M$, then $f$ is locally $L_p$ on $N$
if its pullback $\pi^\ast f$ is locally $L_p$ on $M$.
\end{lemma}

\noindent{\it Proof.}
Choose volume forms $\omega_M$ and $\omega_N$ on $M$ and $N$
respectively. Since $\pi$ is orientation preserving, there is a non-%
negative $C^\infty$ function $\phi(x)$ on $M$ such that
$$
\pi^\ast \omega_N = \phi(x)\omega_M
$$
Now suppose that $x\in N$. Choose a compact neighbourhood $U$ of $x$
in
$N$. Since $\pi$ is proper, $\pi^{-1}(U)$ is compact, and there is
a real number $C$ such that $\phi$ is bounded by $C$ on
$\pi^{-1}(U)$.
Since $\pi^{-1}f$ is locally $L_p$ on $M$, and since $\pi^{-1}(U)$ is
compact,
$$
\int_{\pi^{-1}(U)} |\pi^\ast f|^p \omega_M <\infty.
$$
Consequently,
$$
\int_U |f|^p \omega_N = \int_{\pi^{-1}(U)} |\pi^\ast f|^p
\pi^\ast\omega_N
\le C\int_{\pi^{-1}(U)} |\pi^\ast f|^p \omega_M <\infty. \qed
$$
\medskip

Suppose that $Y$ is a smooth variety which is birational to $X$.
Since
$X$ and $Y$ differ by sets of measure zero, each almost
everywhere defined function on $X$ can be regarded as an almost
everywhere function on $Y$.

\begin{proposition}\label{loc_Lp}
Suppose that $X$ and $Y$ are birational complex algebraic
manifolds. If $p> 0$, then each Bloch-Wigner function on $X$ is
locally $L_p$ when viewed as a function on $Y$.
\end{proposition}

\begin{pf}
If $U$ is a Zariski open subset of $X$, the restriction map
$\bw(X) \to \bw(U)$ is injective. By replacing $X$ by a Zariski open
subset common to $X$ and $Y$, we may assume that $X$ is a Zariski
open subset of $Y$. Let $Z = Y-X$. By Hironaka's resolution of
singularities, there is a smooth variety $\Ytilde$, a normal
crossings divisor $D$ in $\Ytilde$, and a proper map
$$
\pi : (\Ytilde,D) \to (Y,Z)
$$
which induces an isomorphism $\Ytilde - D \to X$. By (\ref{loc_pres})
a function $f\in \bw(X)$ is locally $L_p$ on $Y$ if it is locally
$L_p$
on $\Ytilde$.  It therefore suffices to consider the case where $Z$
is
a normal crossings divisor in $Y$.

Since $X$ is open in $Y$ and Bloch-Wigner functions are smooth on
$X$,
elements of $\bw(X)$ are $L_p$ about points of $X$. Suppose that
$x \in Z$. Choose local coordinates $(z_1,\dots,z_n)$ about $x$
defined
in a polydisk neighbourhood $\Delta$ of $x$ in $Y$
such that $Z$ is contained in the divisor $z_1 z_2 \dots z_n = 0$.

Every multivalued function associated to a relatively closed iterated
integral of holomorphic 1-forms on $\Delta$
with logarithmic singularities along $Z$ an be obtained as follows
(cf. \cite[\S 3]{hain:geom}): There is a $gl_m(\C)$ valued 1-form
$$
\omega \in \Omega^1(\Delta\log Z)\otimes gl_m(\C)
$$
with logarithmic singularities along $Z$ which is integrable:
$$
d\omega + \omega \wedge \omega = 0
$$
and has strictly upper triangular residue along each component of
$Z$.
This defines a flat meromorphic connection on the trivial bundle
$$
\C^m\times \Delta \to \Delta
$$
which is holomorphic on $\Delta - Z$, and has regular singularities
along
$Z$. The multivalued  closed iterated integrals are obtained by
taking linear
combinations of flat sections of this bundle, and then composing with
a linear projection $\C^m \to \C$.

By the several variable generalization of \cite[Theorem 5.5]{wasow}
(see \cite[5.2]{deligne:de}), every flat section of this bundle is of
the form
\begin{equation}\label{standard}
F(z_1,\dots,z_n) =
C P(z_1,\dots,z_n) \prod_{j=1}^n e^{A_j\log z_j},
\end{equation}
where $P: \Delta \to GL_m(\C)$ is a holomorphic function, $C$ is a
constant vector, and $A_j$ is the residue of the connection form
$\omega$ along $z_j = 0$. In particular,  the restriction of an
element of $\tilde{\cal O}(X)$ to an angular sector of $\Delta$ about
0 is an entry of (\ref{standard}). Since the matrices $A_j$ are upper
triangular, it follows that such functions are polynomials in
$\log z_1,\dots, \log z_n$ with coefficients in the ring of
holomorphic functions on $\Delta$. Such functions are $L_p$ on each
closed angular segment of $\Delta$. Elements of $\bw(X)$ are
linear combinations of real and imaginary parts of such functions. It
follows that the restriction of a Bloch-Wigner function to each
closed angular sector of $\Delta$ is also $L_p$ for all $p < 0$. The
result follows.
\end{pf}

Suppose that $f: U^m_n(\C) \to V$ is a continuous function, where $V$
is a finite dimensional real vector space, such as $\R$ or $\R(m)$.
The corresponding function
$$
f^e : GL_m(\C)^{m+n+1} \to V
$$
is an almost everywhere defined cocycle. It should be noted, however,
that
in general this cocycle cannot be made into an everywhere continuous
cocycle. If $f$ is $L_p$ for some $p\ge 1$ when viewed as a function
on
all of the grassmannian, then, by (\ref {loc_Lp}), it will be locally
$L_p$
and locally
integrable with respect to the measure on $GL_m(\C)$ associated to a
left
invariant volume form. Under these conditions, $f$ represents a class
in
$$
\Hp^{m+n}(GL_m(\C),V),
$$
where $\H^\dot$ denotes the locally $L_p$ group cohomology (cf.
\cite{blanc}
and (\ref{coho_def}).)

The natural chain map
$$
j^\ast: C_{\text{cts}}\supdot(GL_m(\C),V)\rightarrow
C\supdot_{\text{loc-}L_p}(GL_m(\C),V)
$$
induces isomorphism between continuous group cohomology and locally
$L_p$ cohomology \cite[(3.5)]{blanc}. We can write down an explicit
inverse of $j^\ast$ as follows. Choose a non-negative continuous
function
$\chi$ on $GL_m(\C)$ with compact support and integral $1$ over the
group.
Define a chain map $r_\chi$ by
\begin{multline*}
r_\chi: C^n_{\text{loc-}L_p}(GL_m(\C),V)\rightarrow
C^n_{\text{cts}}(GL_m(\C),V)\\
\phantom{r9}(r_\chi f)(g_0,\dots,g_n)=\hfill \\
\int_{GL_m(\C)^{n+1}}
\chi(g_0^{-1}h_0)\dots\chi(g_n^{-1}h_n)f(h_0,\dots,h_n)
{\text d}h_0\dots {\text d}h_n.
\end{multline*}
Here $\text{d}h$ denotes a fixed left invariant volume form on
$GL_m(\C)$.
The map $r_\chi$ induces (see \cite[4.11]{blanc}) the inverse
$$
(j^\ast)^{-1}: \Hp\supdot(GL_m(\C),V)
\longrightarrow \Hcts\supdot(GL_m(\C),V)
$$
of $j^\ast$.
In particular, the isomorphism $r_\chi^\ast$ is independent of the
choice
of the bump function $\chi$. Denote by $\mu^\ast$ the composite
$$
\Hp\supdot(GL_m(\C),V)\stackrel{(j^\ast)^{-1}}{\rightarrow}
\Hcts^\ast(GL_m(\C),V)
\rightarrow H^\ast(GL_m(\C),V).
$$
The following result is not unexpected. However, as group cohomology
cycles
have measure zero, and since locally $L_p$ cochains can be changed
with
impunity on sets of measure zero, there is something to prove, and
the result is not immediately obvious.

Recall from (\ref{groupclass}) that a function $f: U^m_n(\C) \to V$
that satisfies the functional equation $A^\ast f = 0$ determines a
class $f^e$ in $H^{m+n}(GL_m(\C),V)$.

\begin{proposition}
Suppose that $f:U^m_n(\C) \to V$ is a continuous function satisfying
the functional equation
$A^\ast f = 0$. If $f$ is locally $L_p$ on the grassmannian for some
$p\ge 1$,
then the class $f^e$ determined by $f$ in $H^{m+n}(GL_m(\C),V)$ is
the
image of the class in $\Hp^{m+n}(GL_m(\C),V)$ determined by $f$ under
the natural map
$$
\mu^{m+n}:\Hp^{m+n}(GL_m(\C);V)\stackrel{\cong}{\leftarrow}
\Hcts^{m+n}(GL_m(\C);V))\rightarrow
H^{m+n}(GL_m(\C);V).
$$
In particular, a generic real Grassmann $m$-logarithm function $D_m$
represents a class in $\Hcts^{2m-1}(GL_m(\C),\R)$
whose image in $H^{2m-1}(GL_m(\C),\R)$ is the class constructed in
(\ref{groupclass}).
\end{proposition}

\begin{pf} From the preceding discussion, we know that the natural
map
$\mu^{m+n}$ is induced by the chain map
$$
\rho_\chi: C^{m+n}_{{\text loc}-L_p}(GL_m(\C), V)
\stackrel{r_\chi}\rightarrow
C^{m+n}_{\text cts}(GL_m(\C),V) \rightarrow C^{m+n}(GL_m(\C),V)
$$
for each bump function $\chi$ on $GL_m(\C)$.

Fix a non-zero vector $e$ in $\C^m$ as a base point. Then $f^e$ is an
$(n+m)-$cocycle in $C^\bullet_{{\text loc}-L_p}(GL_m(\C), V)$. To
distinguish it from $f^e$ viewed as a discrete cocycle, we shall
denote
it by $f^e_{L_p}$. By
(\ref{resoln}), it suffices to prove that the image of $f^e_{L_p}$
under the chain map
$$
C^{m+n}_{{\text loc}-L_p}(GL_m(\C), V)
\stackrel{\rho_\chi}{\rightarrow}
C^{m+n}(GL_m(\C), V) \rightarrow C_e^{m+n}(GL_m(\C), V)
$$
is cohomologous to $f^e$. We prove this by showing that the cocycles
$\rho_\chi f^e_{L_p}$ and $f^e$ agree when evaluated on a cycle with
support
inside $X^{m+n}_{U,e}$.

Choose a sequence of bump functions $\chi_i$ on
$GL_m(\C)$ which converge to the 0-current on $GL_m(\C)$ whose value
on a
test function is its value at the identity.  Since $f^e$ is
continuous at
$x \in X^{m+n}_{U,e}$, then
$$
\lim_{i\rightarrow \infty} \rho_{\chi_i} f^e_{L_p}(x) = f^e(x).
$$
Since value of $\rho_{\chi_i}f^e$ evaluated on a cycle supported on
$X^{m+n}_{U,e}$ is independent of the choice of $\chi_i$, the
proposition follows.
\end{pf}

\section{Chern Classes in Algebraic $K$-theory}
\label{chernclass}

We review the construction of the Chern classes from the
$K$-groups of an affine complex algebraic variety into its Deligne
cohomology. We first fix some notation that will be used throughout
this section.

Let $X$ be an affine variety over $\C$ and $\eta_X$ be its
generic point of $X$. We shall denote the ring of regular functions
$\C[X]$ of $X$ by $R$, and the function field $\C(X)$ of $X$ by $F$.
The
group $GL_m(R)$ will be viewed as the 0-dimensional variety (or, more
accurately, as a direct limit of 0-dimensional varieties) whose
points
are the algebraic functions $f: X \rightarrow GL_m(\C)$. As usual,
$B\subdot GL_m(R)$ will denote the standard simplicial model
of the classifying space of $GL_m(R)$. It will be regarded as a
simplicial
set (or a simplicial variety, each of whose sets of simplices is
0-dimensional).

Let us denote the standard simplicial model of the
universal complex $m$-bundle by
$$
\nu(m): E\subdot GL_m(\C)\rightarrow \bglm.
$$
The pullback of the $l$-th universal Deligne-Beilinson Chern class of
$\nu(n)$
$$
c_l\in \Hd^{2l}(B_\bullet GL_N(\C),\R(l))
$$
along the evaluation map
$$
e_N: X\times B\subdot GL_N(R)\rightarrow B\subdot GL_N(\C)
$$
gives an element $e_N^\ast (c_m)$ of
$$
\Hd^{2m}(X\times B\subdot GL_N(R),\R(m))\cong
\bigoplus_{k=0}^{2l}
\Hd^{2m-k}(X,\R(m))\otimes H^k(GL_N(R),\R).
$$
Evaluating $e_N^\ast(c_m)$ on elements of $H_m(GL_N(R),\R)$
gives a natural map
$$
H_m(GL_N(R),\R)\rightarrow \Hd^{2m-l}(X,\R(m)).
$$
The Chern class
$$
c_{m,l} : K_l(X) \to \Hd^{2m-l}(X,\R(m))
$$
is obtained by taking $n$ to be sufficiently large and composing
the previous map with the Hurewicz homomorphism
$$
K_l(X)  \to H_l(GL_N(R),\R)
$$
Repeating the construction for each Zariski open subset of $X$
and then taking the direct limit over all Zariski open subsets of
$X$,
we obtain the Chern class maps
$$
c_{m,l} : K_l(\eta_X) \rightarrow \Hd^{2m-l}(\eta_X,\R(m)).
$$

\section{Descent of the Universal Deligne-Beilinson Chern
Classes}\label{descent}

In order to prove the existence of the 4-logarithm and all
generic Grassmann logarithms, we prove that the
alternating part, $\Alt c_m$, of the universal Chern class
$$
c_m \in \Hd^{2m}(\bglm, \R(m))
$$
``descends'' to a class
$$
\lambda_m \in s\Hd^{2m}(G^m_\dot, \R(m)).
$$
We begin by making this statement precise. In order
to do this, we need to introduce several simplicial varieties.
In this section, we will distinguish between genuine simplicial
varieties and truncated simplicial varieties. (Cf.\ the conventions
in the introduction.)

First, let
$$
EGL_m(\C) \to BGL_m(\C)
$$
be the standard model of the universal $GL_m(\C)$ bundle in the
category of simplicial varieties. That is, the variety of $n$
simplices,
$E_nGL_m(\C)$, of $EGL_m(\C)$ is $GL_m(\C)^{n+1}$, and the $n$
simplices,
$B_nGL_m(\C)$, of $BGL_m(\C)$ is the quotient of $E_nGL_m(\C)$ by the
diagonal $GL_m(\C)$ action. The face maps of $EGL_m(\C)$ are the
evident
ones. We shall denote this universal bundle by $\nu$.

Define the simplicial variety $E_\dot G^m$ by defining its
variety of $n$-simplices $E_nG^m$ to be
$$
\left\{(v_0,\dots , v_n) : v_j \in \C^m \text{ and the
vectors } v_0, \dots , v_n \text{ are in general position}\right\}.
$$
The $j$th face map is the evident one obtained by forgetting the
$j$th vector. The simplicial variety $B_\dot G^m$ is obtained by
taking
the quotient of $E_\dot G^m$ by the diagonal $GL_m(\C)$ action.
Observe
that its set of $n$ simplices $B_nG^m$ is a point when $n< m$,
and that the projection
$$
E_\dot G^m \to B_\dot G^m
$$
is not a principal $GL_m(\C)$ bundle. Observe also that there is a
natural ``map'' of the $[m,2m]$-truncated simplicial variety
$G^m_\dot$ into $B_\dot G^m$ which is an isomorphism on $n$ simplices
when $m\le n \le 2m$.\footnote{It is tempting, though misleading, to
think
of $G^m_\dot$ as corresponding to a subspace of $B_\dot G^m$.} The
following result is easily proved by considering the double complexes
associated to $G^m_\dot$ and $B_\dot G^m$.

\begin{proposition}\label{identification}
When $m < k \le 2m$, there is a natural isomorphism
$$
\Hd^k(G^m_\dot,\R(m)) \to \Hd^k(B_\dot G^m,\R(m))
$$
which is compatible with the $r \oplus s$ decomposition. \qed
\end{proposition}

Next we want to interpolate between $BGL_m(\C)$ and $B_\dot G^m$.
Fix a non-zero element $e$ of $\C^m$. Define
$$
E_nGL_m(\C)^\gen = \left\{(g_0,\dots,g_n)\in GL_m(\C)^{n+1} :
(g_0e,\dots, g_n e) \in E_nG^m\right\}
$$
and $B_nGL_m(\C)^\gen$ to be the quotient of this by the diagonal
$GL_m(\C)$ action.
We shall denote the principal $GL_m(\C)$ bundle
$$
EGL_m(\C)^\gen \to BGL_m(\C)^\gen.
$$
by $\nu^\gen$. It is the restriction of the universal $GL_m(\C)$
bundle.

By sending $(g_0,\dots,g_n)$ to $(g_0e,\dots,g_ne)$, we obtain a map
$$
\pi : BGL_m(\C)^\gen \to B_\dot G^m.
$$
We have the following diagram of simplicial varieties:
$$
B_\dot G^m \stackrel{\pi}{\leftarrow} BGL_m(\C)^\gen
\hookrightarrow BGL_m(\C)
$$

Denote the component of
$$
c_m(\nu) \in \Hd^{2m}(BGL_m(\C),\R(m))
$$
in
$$
s\Hd^{2m}(BGL_m(\C),\R(m))
$$
by $\Alt c_m(\nu)$. Likewise, we denote the alternating part of
$$
c_m(\nu^\gen)\in \Hd^{2m}(BGL_m(\C)^\gen,\R(m))
$$
by $\Alt c_m(\nu^\gen)$.

We shall identify $\Hd^{2m}(B_\dot G^m,\R(m))$ with
$\Hd^{2m}(G^m_\dot,\R(m))$ via the isomorphism given by
(\ref{identification}). The precise statement of the descent of the
Chern class is:

\begin{theorem}
There is a unique class $\lambda_m$ in $s\Hd^{2m}(G^m_\dot,\R(m))$
such
that $\pi^\ast \lambda_m = \Alt c_m(\nu^\gen)$ in
$s\Hd^{2m}(BGL_m(\C)^\gen,\R(m))$.
\end{theorem}

\begin{remark}
The theorem and our proof are equally valid with $\Q(m)$
coefficients.
\end{remark}

The remainder of this section is devoted to the proof of this
theorem.
Because the projection $E_\dot G^m \to B_\dot G^m$ is not a principal
bundle, it is convenient to introduce a simplicial variety which
interpolates
between $B_\dot G^m$ and $BGL_m(\C)^\gen$. Define $\Btilde_\dot G^m$
to
be the homotopy quotient
$$
\left(EGL_m(\C) \times E_\dot G^m \right)/GL_m(\C)
$$
of $E_\dot G^m$ by $GL_m(\C)$. Set
$$
\Etilde_\dot G^m = EGL_m(\C) \times E_\dot G^m.
$$
Then the natural projection $\Etilde_\dot G^m \to \Btilde_\dot G^m$
is a principal $GL_m(\C)$ bundle which we shall denote by $\mu$. The
projection
$$
EGL_m(\C) \times E_\dot G^m \to E_\dot G^m
$$
induces a morphism $p : \Btilde_\dot G^m \to B_\dot G^m$.

Set
$$
\Etilde GL_m(\C)^\gen = EGL_m(\C) \times EGL_m(\C)^\gen
$$
and $\Btilde GL_m(\C)^\gen$ equal to the quotient of this by the
diagonal $GL_m(\C)$ action.
We have the diagram
$$
\begin{matrix}
EGL_m(\C)^\gen & \leftarrow & \Etilde GL_m(\C)^\gen &
\to & \Etilde_\dot G^m \cr
\downarrow  & & \downarrow & & \downarrow \cr
BGL_m(\C)^\gen & \leftarrow & \Btilde GL_m(\C)^\gen & \to &
\Btilde_\dot G^m & \to & B_\dot G^m \cr
\end{matrix}
$$
of simplicial varieties where the vertical arrows are principal
$GL_m(\C)$ bundles, and where the right hand map in the top row is
induced by evaluation on $e \in \C^m - \{0\}$.
The morphism $\Btilde GL_m(\C)^\gen \to BGL_m(\C)^\gen$ is a homotopy
equivalence of simplicial varieties, and therefore induces an
isomorphism
on Deligne cohomology. The class $c_m(\nu^\gen)$ therefore descends
naturally to the class
$$
c_m(\mu) \in \Hd^{2m}(\Btilde_\dot G^m,\R(m)).
$$
We will prove the theorem by showing that there is a class
$$
\lambda_m \in  s\Hd^{2m}(B_\dot G^m,\R(m))\cong
s\Hd^{2m}(G^m_\dot,\R(m))
$$
such that
$$
p^\ast \lambda_m = \Alt c_m(\mu) \in s\Hd^{2m}(\Btilde_\dot
G^m,\R(m)).
$$

Observe that each of the varieties defined in this section is a
$\Sigma_\dot$ variety and that all morphisms between them that we
have constructed in this section respect the $\Sigma_\dot$
structures.

Denote the $(m-1)$ skeleton of $\Btilde_\dot G^m$ by $\Btilde_{<m}
G^m$.
(This is the $[0,m-1]$ truncated simplicial variety whose $n$
simplices
are  identical with those of $\Btilde_\dot G^m$ when $n < m$ and
empty
otherwise.)

\begin{proposition}\label{les}
There is a long exact sequence
\begin{multline*}
\dots \to \Hd^{k}(G^m_\dot,\R(m)) \to \Hd^{k}(\Btilde_\dot G^m,\R(m))
\to \Hd^{k}(\Btilde_{<m}G^m,\R(m)) \to \\
\dots \to \Hd^{2m}(G^m_\dot,\R(m)) \to \Hd^{2m}(\Btilde_\dot
G^m,\R(m))
\to \Hd^{2m}(\Btilde_{<m}G^m,\R(m))
\end{multline*}
which remains exact when the alternating part functor $s$ is applied.
\end{proposition}

\begin{pf}
Let $B_{\ge m}G^m$ be the $[m,\infty)$-truncated simplicial variety
whose $n$ simplices are those of $B_\dot G^m$ when $n\ge m$ and empty
otherwise. Let
$\Btilde_{\ge m} G^m$ be the analogous $[m,\infty)$-truncated
simplicial 
variety constructed out of the $n$ simplices of $\Btilde_\dot G^m$
for
$n\ge m$. The natural projection $\Btilde_{\ge m} G^m \to B_{\ge
m}G^m$
induces an isomorphism on Deligne cohomology as $\Btilde_n G^m \to
B_n G^m$ is a homotopy equivalence whenever $n \ge m$.

An easy spectral sequence argument shows that the inclusion
$G^m_\dot \to B_{\ge m}G^m$ induces an isomorphism
$$
\Hd^k(B_{\ge m}G^m,\R(m)) \to \Hd^{k}(G^m_\dot,\R(m))
$$
when $k\le 2m$. This isomorphism is compatible with symmetric group
actions. We therefore have an isomorphism
$$
\Hd^{k}(G^m_\dot,\R(m)) \cong \Hd^k(\Btilde_{\ge m}G^m,\R(m))
$$
when $k\le 2m$, also compatible with symmetric group actions.

Finally, observe that the sequence
$$
0 \to \R_\calD^\dot(\Btilde_{\ge m}G^m) \to
\R_\calD^\dot(\Btilde_\dot G^m)
\to \R_\calD^\dot(\Btilde_{<m}G^m) \to 0
$$
of Deligne cochain complexes is exact and compatible with the
symmetric
group actions. It induces a long exact sequence on cohomology.
The result follows by combining these results.
\end{pf}

Since the sequence
$$
s\Hd^{2m}(G^m_\dot,\R(m)) \to s\Hd^{2m}(\Btilde_\dot G^m,\R(m))
\to s\Hd^{2m}(\Btilde_{<m}G^m,\R(m))
$$
is exact, the existence of a lift $\lambda_m$ of $\Alt c_m(\mu^\gen)$
will be proved if we can show that the image
of $\Alt c_m(\mu)$
in
$$
s\Hd^{2m}(\Btilde_{<m}G^m,\R(m))
$$
vanishes.

\begin{proposition}\label{triv_act}
If $0 \le n < m$, then $\Hd^\dot(\Btilde_nG^m,\R(m))$ is a trivial
$\Sigma_{n+1}$-module.
\end{proposition}

\begin{pf} We begin the proof with an elementary observation. Suppose
that
$A$ is an $m\times k$ matrix of complex numbers.  If $k\le m$ and if
the
columns of $A$ are linearly independent, then so are the columns of
$AB$
for all $B\in GL_k(\C)$. This is not the case when $k>m$: if each $m$
columns of $A$ are linearly independent, then it is not true that
each $m$
of the columns of $AB$ are linearly independent for all $B\in
GL_k(\C)$.

In other notation, this says that $E_nG^m$ has a natural right
action of $GL_{n+1}(\C)$ provided that $n<m$. This action commutes
with the diagonal left action of $GL_m(\C)$ on $E_nG^m$. After
taking the product with $EGL_m(\C)$ and taking the quotient by
$GL_m(\C)$, we see that $\Btilde_nG^m$ has a natural right
$GL_{n+1}(\C)$ action, provided $n< m$.

Identify $\Sigma_{n+1}$ with the subgroup of $GL_{n+1}(\C)$
consisting
of all permutation matrices. The restriction of the right action
of  $GL_{n+1}$ on $E_nG^m$ to $\Sigma_{n+1}$ is its standard action.
It follows that the action of $\Sigma_{n+1}$ on $\Btilde_nG^m$ is the
restriction of the $GL_{n+1}(\C)$ action. Since $GL_{n+1}(\C)$ is
connected, it follows that the automorphism of $\Btilde_nG^m$
induced by an element of $\Sigma_{n+1}$ is homotopic to the
identity.
\end{pf}

\begin{proposition}\label{equivalence}
If $n<m$, then there is a natural map of simplicial varieties
$$
BGL_{m-n-1}(\C) \to \Btilde_nG^m
$$
which is a homotopy equivalence.
\end{proposition}

\begin{pf}
This follows from two facts: First, $GL_m(\C)$ acts
transitively on $E_nG^m$ and the isotropy group of a point is
$$
G(n) = {\left (
\begin{array}{cc}
 I_{n+1} & \ast \\
 0 & {\rm GL}_{m-n-1}(\C)
\end{array} \right )}.
$$
Second, the inclusion of $GL_{m-n-1}(\C)$ into $G(n)$ is a
homotopy equivalence.
\end{pf}

\begin{lemma}\label{isomorphism}
The inclusion $\Btilde_0G^m \hookrightarrow \Btilde_{<m}G^m$
induces an isomorphism
$$
s\Hd^\dot(\Btilde_{<m}G^m, \R(m)) \cong \Hd^\dot(\Btilde_0G^m,\R(m)).
$$
Consequently, there is a natural isomorphism
$$
s\Hd^\dot(\Btilde_{<m}G^m,\R(m))\cong \Hd^\dot(BGL_{m-1}(\C),\R(m)).
$$
\end{lemma}

\begin{pf} The first isomorphism follows immediately from
(\ref{triv_act}) by looking at the spectral sequence associated to
$\Btilde_{<m}G^m$. The second assertion follows the previous result.
\end{pf}

To complete the proof of the theorem, observe that the restriction of
$\Etilde G^m \to \Btilde G^m$ to $\Btilde_0 G^m$ has structure group
the group $G(0)$ defined in the proof of (\ref{equivalence}). Since
this
group is homotopy equivalent to $GL_{m-1}(\C)$, it follows that the
image of $c_m(\mu)$ in
$$
\Hd^{2m}(\Btilde_0G^m,\R(m)) = s\Hd^{2m}(\Btilde_{<m},\R(m))
$$
vanishes. This establishes the existence of $\lambda_m$.

To prove uniqueness, note that it follows from (\ref{les}) that
$\lambda_m$ is unique if
\begin{equation}\label{silly}
s\Hd^{2m-1}(\Btilde_\dot G^m,\R(m)) \to
s\Hd^{2m-1}(\Btilde_{<m}G^m,\R(m)).
\end{equation}
is surjective. By (\ref{isomorphism}),
$$
s\Hd^{2m-1}(\Btilde_{<m}G^m,\R(m)) \cong
H^{2m-2}(BGL_{m-1}(\C),\C/\R(m)).
$$
Thus, to prove that (\ref{silly}) is surjective, it suffices to prove
that the restriction mapping
$$
H^{2m-2}(\Btilde_\dot G^m) \to H^{2m-2}(\Btilde_0 G^m)\cong
H^{2m-2}(BGL_{m-1}(\C)).
$$
induced by the inclusion $\Btilde_0 G^m \hookrightarrow \Btilde_\dot
G^m$
is surjective. This follows as the restriction of the natural
$GL_m(\C)$
bundle $\mu$ to $\Btilde_0 G^m$ corresponds to the universal
$GL_{m-1}(\C)$ bundle over $BGL_{m-1}(\C)$. The Chern classes
$c_1(\mu),\dots,
c_{m-1}(\mu)$ therefore restrict to the generators of the cohomology
ring
of $\Btilde_0 G^m$. Surjectivity follows and, along with it, the
uniqueness
of $\lambda_m$.

\section{Chern Classes in Algebraic $K$-theory---Addendum}

In this section we prove two results needed in the proof of the
existence of Grassmann logarithms and in relating them to Chern
classes
on algebraic $K$-theory. The first result asserts that the
class
$$
\lambda_m \in s\Hd^{2m}(G^m_\dot,\R(m))
$$
can be used to represent the restriction
$$
c_m : r_mK_p(\eta_X) \to \Hd^{2m-p}(\eta_X,\R(m))
$$
of the Chern class to the rank $m$ part of the algebraic $K$-theory
of
the generic point $\eta_X$ of each complex algebraic variety $X$. The
second result asserts that the restriction of the class $\lambda_m$
to
$G^m_0$ is the volume form $\vol_m$.

Let $U$ be a smooth Zariski open subset of $X$. Denote its
coordinate ring by $\C[U]$. Let $G^m_\dot(\C[U])$ denote the
simplicial
set (i.e., 0-dimensional simplicial variety) whose $n$ simplices
consist
of all regular maps $U \to G^m_n$. The evaluation map
$$
U \times G^m_\dot(\C[U]) \to G^m_\dot
$$
induces a map
$$
\Hd^{2m}(\bgm,\R(m)) \to \Hd^{2m}(U\times \bgm(\C[U]),\R(m))
$$
on Deligne cohomology. As in the case of the construction of Chern
classes
on $K$-theory, by evaluation on $\lambda_m$ we obtain maps
$$
l_{p,m} : H_p(GL_m(\C[U])) \to \Hd^{2m-p}(U,\R(m)).
$$
Denote the function field of $X$ by $F$.

\begin{theorem}\label{lambda-chern}
The maps $l_{p,m}$ induce the restriction of the $m$-th Chern class
$$
c_m : r_mK_p(\eta_X) \to \Hd^{2m-p}(\eta_X,\R(m))
$$
to the rank $m$ part of $K_p(\eta_X)$. That is, if $x\in
r_mK_p(\eta_X)$,
then
$$
c_m(x) = l_{p,m}(\tilde{x})
$$
where $\tilde{x}$ is any element of $H_p(GL_m(F),\Q)$ whose image in
$$
K_p(\eta_X)_\Q \subseteq H_p(GL(F),\Q)
$$
is $x$.
\end{theorem}

\begin{remark}
This construction  (and the theorem) are equally valid for the class
$\lambda_m|_{U^m_\dot}$ in $\Hd^{2m}(U^m_\dot,\R(m))$,
 where $U^m_\dot$ is a Zariski open subset of $G^m_\dot$ that
satisfies
the condition (\ref{condition}) of \S \ref{polylog} --- see
(\ref{variant}).
\end{remark}

\begin{pf} We begin by showing that elements of
$\Hd^{2m}(BGL_m(\C)^\gen,\R(m))$ also induce maps
$$
H_p(GL_m(\C[U])) \to \Hd^{2m-p}(U,\R(m))
$$
for all smooth varieties. The construction is very similar to that of
the
universal Chern classes and the maps $l_{p,m}$, so we'll be brief.

View $BGL_m(\C[U])^\gen$ as the simplicial set whose $n$ simplices
consist of all regular maps
$$
U \to B_nGL_m(\C)^\gen.
$$
The classes $\Alt c_m(\nu^\gen)$ and $c_m(\nu_\gen)$ both induce maps
$$
H_p(GL_m(\C[U])) \to \Hd^{2m-p}(U,\R(m)).
$$
It follows immediately from (\ref{sign-copy}) that these two maps
agree.
Denote this map by $c_m^\gen$.

By the naturality of the constructions, the diagram (whose horizontal
maps are induced by evaluation)
$$
\begin{CD}
\Hd^{2m}(\bgm,\R(m)) @>>> \Hd^{2m}(U\times \bgm(\C[U]), \R(m))\\
@VVV                      @VVV \\
\Hd^{2m}(BGL_m(\C)^{gen},\R(m)) @>>> \Hd^{2m}(U\times
BGL_m(\C[U])^{gen},\R(m))\\
@AAA                      @AAA \\
\Hd^{2m}(BGL_M(\C),\R(m)) @>>> \Hd^{2m}(U\times BGL_M(\C[U]),\R(m))
\end{CD}
$$
commutes for all $M\ge m$. By taking $M$ to be sufficiently large
($M\ge p$ will do by Suslin \cite{suslin}), and taking the limit
over all smooth open subsets $U$ of $X$, we see that the diagram
$$
\begin{CD}
H_p(GL_m(F),\Q) @>l_{p,m}>> \Hd^{2m-p}(\eta_X,\R(m))\\
@VVV                      @VVV \\
H_p(GL_m(F),\Q) @>c_m^\gen>> \Hd^{2m-p}(\eta_X,\R(m))\\
@AAA			  @AAA \\
H_p(GL(F),\Q)   @>c_m>>  \Hd^{2m-p}(\eta_X,\R(m))\\
@A\text{Hurewicz}AA        @AAA\\
K_p(\eta_X)_\Q  @>c_m>>  \Hd^{2m-p}(\eta_X,\R(m))
\end{CD}
$$
commutes. The result follows.
\end{pf}

\begin{remark}\label{factorization}
The proof actually shows that $\lambda_m$ induces a map
$$
\lbar_{p,m} : GH^m_p(\eta_X) \to \Hd^{2m-p}(\eta_X,\R(m))
$$
and that $l_{p,m}$ is the composition with $\lbar_{p,m}$ with
the natural map
$$
H_p(GL_m(F)) \to GH^m_p(\eta_X).
$$
\end{remark}

Next, we determine the restriction of the class $\lambda_m$ to
$G^m_0$.

\begin{theorem}
\label{multiple}
The image of $\lambda_m$ under the restriction mapping
$$
s\Hd^{2m}(G^m_\bullet,\R(m))\rightarrow sH^m(G^m_0,\C)
$$
is $\vol_m$.
\end{theorem}

Suppose that $k$ is a field. In the rest of this section, we shall
denote the $K$-theory and Grassmann homology of $\spec k$ by
$K_\dot(k)$ and $GH^m_\dot(k)$, respectively.

Before giving the proof, we review some results of Suslin from
\cite{suslin}. For this discussion, $k$ is an infinite field.  Define
$S_m(k)$ to be
$$
\left[\coker\Bigg\{ \bigoplus_{E_{m+2}G^m(k)} \Z
\stackrel{A_\ast}{\longrightarrow}
\bigoplus_{E_{m+1}G^m(k)}\Z\Bigg\}\right]\otimes_{GL_m(k)}\Z.
$$
This is the group of Grassmann $m$-chains mod boundaries. It is
generated by the equivalence class of the $(m+1)$ tuples of vectors
$$
(e_1,\dots, e_m, \sum a_ie_i),
$$
where $e_1,\dots, e_m$ is the standard basis of $k^m$. Following
Suslin, we shall denote the corresponding element of $S_m(k)$ by
$\langle a_1,\dots, a_m \rangle$.
There is a natural inclusion
\begin{equation}\label{into_S}
GH^m_m(k) \hookrightarrow S_m(k)
\end{equation}
whose cokernel is trivial or $\Z$ according to whether $m$ is odd or
even. By (\ref{grass_homo}), there is a map
\begin{equation}\label{into_GH}
H_m(GL_m(k)) \to GH^m_m(k).
\end{equation}

Denote the Milnor $K$-theory of $k$ by $K^M_\dot(k)$. It is a ring
generated by $K^M_1(k) = k^\times$. The symbol
$\{a_1,\dots,a_m\}\in K^M_m(k)$ is the product of the
$a_i\in k^\times$.

Suslin shows that there is a well defined homomorphism
\begin{equation}\label{into_milnor}
S_m(k) \to K^M_m(k)
\end{equation}
defined by
$$
\langle a_1,\dots, a_m \rangle \mapsto \{a_1,\dots,a_m\}.
$$

Suslin also proves that the map
$$
H_m(GL_m(k)) \to H_m(GL(k))
$$
is an isomorphism. Consequently, the Hurewicz homomorphism
indues a homomorphism
$$
K_m(k) \to H_m(GL_m(k))
$$
Composing this with (\ref{into_S}), (\ref{into_GH}), and
(\ref{into_milnor}), he obtains a map
$$
\phi : K_m(k) \to K^M_m(k).
$$
He proves that it has the property that if $a_1,\dots,a_m \in
k^\times$, then
$$
\phi(a_1\cdot \dots \cdot a_m) =
(-1)^{m-1}(m-1)! \{a_1,\dots,a_m\}.
$$

\begin{pf*}{Proof of Theorem \ref{multiple}}
There is a rational number $K$ such that the image of $\lambda_m$
in $H^m(G^m_0,\Q(m))$ is $K\vol_m$. Let $X = \C^m$. We will
determine $K$ by evaluating $l_{m,m}$ on a class in $K_m(\eta_X)$
and comparing the answer with the value of the Chern class on it.

Deligne cohomology maps to de~Rham cohomology, and it will be
sufficient for our needs to replace the Chern class and the map
$l_{m,m}$ with their composite with the map to de~Rham
cohomology.

Denote the function field of $X$ by $F$. The de~Rham cohomology of
$\spec F$ is the set of K\"ahler differentials $\Omega_{F/\C}^p$. It
is
a standard fact that the first Chern class
$$
c_1 : K_1(F) \to \Omega^1_{F/\C}
$$
takes $f$ to $df/f$. It follows from standard properties of Chern
classes on algebraic $K$-theory (see, for example,
\cite[p.~28]{schneider}), that
\begin{equation}\label{chern_map}
c_m : K_m(F) \to \Omega^m_{F/\C}
\end{equation}
takes $f_1\cdot \dots \cdot f_m$ ($f_i \in F^\times$) to
$$
(-1)^{m-1}(m-1)! \frac{df_1}{f_1}\wedge \dots \wedge
\frac{d f_m}{f_m}
$$
It follows that if we define
$$
\psi : K^M_m(F) \to \Omega^m_{F/\C}
$$
by
$$
\psi : \{f_1,\dots ,f_m\} \to
\frac{df_1}{f_1}\wedge \dots \wedge\frac{d f_m}{f_m},
$$
then the restriction
$$
K^M_m(F) \to K_m(F) \stackrel{c_m}{\to} \Omega^m_{F/\C}
$$
of the Chern class to the Milnor $K$-theory is the composite
$$
K^M_m(F) \to K_m(F) \stackrel{\phi}{\to} K^M_m(k)
\stackrel{\psi}{\to} \Omega^m_{F/\C}.
$$

To compute the constant $K$ that relates $\lambda_m$ and $\vol_m$,
we use the fact that the map
$$
H_m(GL_m(F)) \to H_m(GL(F))
$$
is an isomorphism. From (\ref{lambda-chern}), it follows that we can
compute (\ref{chern_map}) using $\lambda_m$.

Observe that since $S_m(F)$ is the set of all Grassmann $m$-chains
mod boundaries, the class $\lambda_m$ induces a map
$$
S_m(F) \to \Omega^m_{F/\C}
$$
such that the diagram
$$
\begin{CD}
K_m(F) @>>> S_m(F) @>\lambda_m>> \Omega^m_{F/\C} \\
@V\text{Hurewicz}VV              @AA\text{inclusion}A \\
H_m(GL_m(F)) @>>> GH^m_m(F) \\
\end{CD}
$$
commutes.

To compute the map $S_m(F) \to \Omega^m_{F/\C}$
induced by $\lambda_m$, it is first necessary to realize that there
are two descriptions of $G^m_0$: first, it is the quotient of the set
of
$(m+1)$-tuples of vectors  $(v_0,v_1,\dots, v_m)$ in $\C^m$,  in
general position, mod the diagonal action $GL_m$. In this case an
isomorphism with $(k^\times)^m$ is given by
$$
(a_1,\dots,a_m) \mapsto (e_1,\dots,e_m,\sum a_i e_i).
$$
The second description is $\P^m$ minus the union of the coordinate
hyperplanes. This is identified with $(\C^\times)^m$ via the
formula
$$
(x_1,\dots,x_m) \mapsto [1,x_1,\dots,x_m].
$$
The two descriptions are related by identifying orbit of the vectors
$(v_0,v_1,\dots, v_m)$ with the point of $\P^m$ corresponding to the
kernel of the linear map $k^{m+1} \to k$ that takes $e_i$ to $v_i$
---
cf.\ \cite[\S 5]{hain-macp}.

The volume form $\vol_m$ on $G^m_0$ is, by convention (cf.
\cite[p.~422]{hain-macp}),
$$
\frac{dx_1}{x_1} \wedge \dots \wedge \frac{dx_m}{x_m}.
$$
A short computation then shows that this equals the form
$$
\frac{da_1}{a_1} \wedge \dots \wedge \frac{da_m}{a_m}
$$
with respect to the quotient coordinates. It follows that the
map
$$
S_m(F) \to \Omega^m_{F/\C}
$$
induced by $\lambda_m$ takes $\langle a_1,\dots,a_m \rangle$ to
$$
K\frac{da_1}{a_1} \wedge \dots \wedge \frac{da_m}{a_m}.
$$
Since $l_{m,m}$ is the composite
$$
H_m(GL_m(F)) \to GH_m(F^m) \to H^m(\eta_X,\R(m)),
$$
and since $l_{m,m}$ equals $c_m$, we deduce that $K= 1$.
\end{pf*}

\section{Generic Grassmann Polylogarithms---Existence and Relation
to Chern Classes}

We first use the results of the preceding sections to prove the
existence of generic real Grassmann logarithms and establish their
relation to the Beilinson Chern classes.

\begin{theorem}
For all $m$, there is a canonical choice of a generic real Grassmann
$m$-logarithm. Moreover, for all complex algebraic varieties $X$,
this
$m$-logarithm induces the restriction
$$
c_m : r_mK_n(\eta_X) \to \Hd^{2m-n}(\eta_X,\R(m))
$$
of the Beilinson-Chern class to the $m$th part of the rank filtration
of $K_\dot(\eta_X)$.
\end{theorem}

\begin{pf}
By \cite[(7.1)]{hain:generic}, there is a Zariski open subset
$V^m_\dot$ of $G^m_\dot$ where $U^m_0 = G^m_0$ and where each
$V^m_n$ is a rational $K(\pi,1)$. By (\ref{summand}), there is a
canonical injection
$$
\Hd^{2m}(V^m_\dot,\R(m)) \hookrightarrow
H_\bw^{2m}(V^m_\dot,\R(m)).
$$
Let $L_m'$ be the restriction of $\lambda_m$ to $V^m_\dot$, viewed
as a class in $H_\bw^{2m}(V^m_\dot,\R(m))$. We need to skew
symmetrize. Let $U^m_\dot$ be the Zariski open subset where
$U^m_n$ is the intersection of the translates of $V^m_n$ under the
action of $\Sigma_{m+n+1}$ on $G^m_n$. Then $U^m_0=G^m_0$ and
$U^m_\dot$ is a $\Sigma_\dot$ variety. Let
$$
L_m \in sH_\bw^{2m}(U^m_\dot,\R(m))
$$
be the alternating part of the restriction of $L_m'$ to $U^m_\dot$.
It
follows from (\ref{multiple}) that $L_m$ is a generic Grassmann
$m$-logarithm.

The second assertion is an immediate consequence of the definition
of Grassmann logarithms, Theorems \ref{lambda-chern} and
\ref{multiple}, and the fact that the open subvarieties $U^m_\dot$ of
$G^m_\dot$ used above always satisfy the condition (\ref{condition})
of \S \ref{polylog}.
\end{pf}

\section{Comparison of Cohomologies}
\label{comparison}

In this section we prove Theorem \ref{summand}.
We begin with a brief guide to the proof. For each
smooth variety $X$ with $q=0$, we will construct a functorial
complex $\R(m)^\dot_\F(X)$ which is a formal analogue of the
complex $\R(m)^\dot_\bw(X)$. There will be a natural chain map
$$
\R(m)^\dot_\F(X) \to \R(m)^\dot_\bw(X).
$$
The homology of the formal $\bw$-complex will be denoted by
$\Hf^\dot(X,\R(m))$.  Taking homology, we will have a commutative
diagram
$$
\begin{matrix}
\Hf^\dot(X,\R(m)) & \to & H_\bw^\dot(X,\R(m)) \cr
&\searrow & \downarrow \cr
& & \Hd^\dot(X,\R(m)) \cr
\end{matrix}
$$
We will show, when $X$ is an rational $n$-$K(\pi,1)$, that the
composite
$$
\Hf^\dot(X,\R(m)) \to \Hd^\dot(X,\R(m))
$$
is an isomorphism in degrees $\le n$ and injective in dimension
$n+1$.
The result in the case when
$X_\dot$ is a single space then follows. The simplicial version will
follow using a spectral sequence argument.

Our first task is to construct the complex $\R(m)_\F^\dot(X)$. To do
this, we need to construct a formal analogue $\Omega^\dot_\R(X)_\F$
of
$\Omega_\R^\dot(X)$ and a formal analogue $\bw(X)_\F$ of the
ring of $\bw(X)$.

Since $q(X)=0$, it follows from elementary Hodge theory that there
are regular functions $f_j : X \to \C$ such that $\Omega^1(X)$ has
basis $df_1/f_1,\dots, df_m/f_m$.

Let $A_\C^\dot(X)$ be the $\C$-subalgebra of $\Omega^\dot(X)$
generated by the $df_j/f_j$.  Let
$$
\theta_j = d\Arg f_j \text{ and } \rho_j = d \log|f_j|.
$$
Note that $df_j/f_j = \rho_j + i \theta_j$. Observe that
$\Lambda^\dot_\C(df_j/f_j:j = 1,\dots, m)$  is a subalgebra of
$\Lambda_\C^\dot(\theta_j,\rho_j: j=1,\dots,m)$. The latter algebra
has
the real form $\Lambda_\R^\dot(\theta_j,\rho_j)$. Each element $u$
of the ideal
$$
K :=
\ker\left\{ \Lambda_\C^\dot (df_j/f_j:j=1,\dots,m) \to
A^\dot(X)\right\}
$$
can be viewed as elements of
$\Lambda_\C^\dot(\theta_j,\rho_j:j=1,\dots,m)$.
So we can write each such $u$ in the form $a(u) + i b(u)$, where
$a(u), b(u) \in \Lambda_\R^\dot(\theta_j,\rho_j: j = 1,\dots,m)$.
We define $\Omega_\R^\dot(X)_\F$ to be the algebra
$$
\Lambda_\R^\dot(\theta_j,\rho_j: j=1,\dots, m)/(a(u), b(u): u\in K).
$$

The ring of formal Bloch-Wigner functions $\bw(X)_\F$ is defined in
terms of the Malcev completion of the fundamental group of $X$.
Denote the complex form of the Malcev group associated to
$\pi_1(X,x)$
by $G(X,x)$. Denote its real form by $G_\R(X,x)$.  Each of these is
the inverse limit of its finite dimensional quotients $G(X,x)_s$.
The quotient $G_\R(X,x) \backslash G(X,x)$ is a real proalgebraic
variety. Its  coordinate ring is, by definition, the direct limit of
the
coordinate rings of its canonical quotients
$G_\R(X,x)_s \backslash G(X,x)_s$

Recall that each path $\gamma$ in $X$ from $x$ to $y$ induces a
group isomorphism $\mu_\gamma : G(X,x) \to G(X,y)$ which
preserves real forms.

\begin{proposition}\label{translation}
(a) For all $x\in X$, there is a canonical real analytic map
$$
\mu_x : X \to G_\R(X,x) \backslash G(X,x)
$$
(b) If $\gamma$ is a path in $X$ from $x$ to $y$, then the
diagram
$$
\begin{matrix}
X &\stackrel{\mu_x}{\longrightarrow} & G_\R(X,x) \backslash
G(X,x)\cr
& \mu_y \searrow & \downarrow \mu_\gamma \cr
&&  G_\R(X,y) \backslash G(X,y)\cr
\end{matrix}
$$
commutes, where the vertical map is the one induced by
$\mu_\gamma$.  \newline
(c) If $\gamma$ is a loop in $X$ based at $x$, then
$$
\mu_\gamma^\ast : \O_\R\left(G_\R(X,x) \backslash G(X,x)\right) \to
\O_\R\left(G_\R(X,x) \backslash G(X,x)\right),
$$
is the identity. Here $\O_\R(Y)$ denotes the coordinate ring of
the real proalgebraic variety $Y$.
\end{proposition}

\begin{pf}
We will use the notation and terminology of \cite[\S 7]{hain-macp}.
Observe that $G_\R(X,x)_s$ is the real Zariski closure of $\Gamma_s$
in $G_s$. The map $\mu_x$ is simply the inverse limit of the
composites
$$
X \stackrel{\theta_x^s}{\longrightarrow} \Alb^s(X,x) \to
G_\R(X,x)_s\backslash G_s.
$$
If $\gamma$ is a path in $X$ from $x$ to $y$, then we have the
element
$T_s(\gamma)$ of $G_s$. The sequence $\left\{T_s(\gamma)\right\}$
converges to an element $T(\gamma)$ of $G$. We have
$$
G_\R(X,y) = T(\gamma)^{-1}G_\R(X,x)T(\gamma).
$$
The map $\mu(\gamma)$ is induced by left multiplication by
$T(\gamma)$.
The final statement follows as the coordinate ring of
$G_\R(X,x) \backslash G(X,x)$ is the ring of functions on $G(X,x)$
that are invariant under left multiplication by elements of
$G_\R(X,x)$.
\end{pf}

Define the ring of {\it formal Bloch-Wigner functions on $X$},
$\bw(X)_\F$, to be the coordinate ring of $G_\R(X,x) \backslash
G(X,x)$. It follows from (\ref{translation}) that the assignment of
$\bw(X)_\F$ to $X$ is a well defined contravariant functor.

\begin{remark}
What is called the ring of Bloch-Wigner functions in
\cite[\S 11]{hain-macp} is what we are defining to be the ring of
formal Bloch-Wigner functions in this paper.
\end{remark}\medskip

The {\it formal Bloch-Wigner complex of $X$} is defined by
$$
\Omega^\dot_\bw(X)_\F = \bw(X)_\F \otimes_\R \Omega^\dot_\R(X)_\F.
$$
It is a differential graded $\R$-algebra canonically associated to
$X$.

There are natural weight filtrations on $A^\dot_\C(X)$,
$\Omega_\R^\dot(X)_\F$ and $\bw(X)_\F$, and therefore on
$\Omega^\dot_\bw(X)_\F$. These are defined as follows: The weight
filtration on $A^\dot_\C(X)$ is the one induced by its inclusion into
$\Omega^\dot(X)$. Since each $df_j/f_j$ has
weight 2, it follows that all elements of $A^m_\C(X)$ have weight 2m.
The weight filtration on $\Omega_\R^\dot(X)_\F$ is also defined this
way --- all elements of degree $m$ have weight $2m$. The weight
filtration on $\bw(X)_\F$ is defined in \cite[\S 11]{hain-macp}. It
is not difficult to check that the weight filtration of
$\Omega^\dot_\bw(X)_\F$ is a filtration by subcomplexes. (Use
\cite[(7.7)]{hain-macp}.)

Finally, we define $\Hf^\dot(X,\R(m))$, the {\it formal Bloch-Wigner
cohomology of $X$ with coefficients in $\R(m)$}, to be the cohomology
of the complex
$$
\R(m)_\F^\dot(X) = \cone[F^pW_{2p}A_\C\supdot(X)
\rightarrow W_{2p}\obw\supdot(X)_\F\otimes\R(m-1)][-1].
$$
The map $A_\C^\dot(X) \to \obw^\dot(X)_\F\otimes\R(m-1)$ is the
composite of
$$
A^\dot_\C(X) \to \Omega^\dot_\R(X)_\F\otimes\left(\R(m-1) \oplus
\R(m)\right) \to \Omega^\dot_\R(X)_\F\otimes\R(m-1),
$$
where the first map is the algebra homomorphism that takes $df_j/f_j$
to $\rho_j + i \theta_j$, with the natural inclusion
$$
\Omega^\dot_\R(X)_\F\otimes\R(m-1) \hookrightarrow
\obw^\dot(X)_\F\otimes\R(m-1).
$$

When $X_\dot$ is a simplicial complex algebraic manifold where each
$X_m$
has $q=0$, we define $\R(m)_\F^\dot(X_\dot)$ to be the total complex
associated to the cosimplicial chain complex obtained by applying the
functor $\R(m)_\F^\dot(\blank)$ to $X_\dot$.

It follows from \cite[pp.~436--7]{hain-macp} that there is a natural
$W_\dot$ filtered algebra homomorphism $\bw(X)_\F \to \bw(X)$.
There is an obvious filtered algebra homomorphism
$\Omega^\dot_\R(X)_\F \to \Omega^\dot_\R(X)$. These induce chain maps
$$
\R(m)_\F^\dot(X) \to \R(m)_\bw^\dot(X) \to \dbcomplex^\dot(X).
$$
Similarly, in the simplicial case, we have chain maps
$$
\R(m)_\F^\dot(X_\dot) \to \R(m)_\bw^\dot(X_\dot) \to
\dbcomplex^\dot(X_\dot).
$$

Theorem \ref{summand} will follow directly from the following result.

\begin{theorem}\label{tech}
Suppose that $X_\dot$ is a simplicial complex algebraic manifold
where
each $X_m$ has $q=0$. If, for all $m$, $X_m$ is a rational
$(n-m)$-$K(\pi,1)$, then the natural map
$$
\Hf^t(X_\dot,\R(m)) \to \Hd^t(X_\dot,\R(m))
$$
is an isomorphism when $t \le n$.
\end{theorem}

The proof of this result occupies the rest of this section.  The
first step is to observe that it follows from the analogue of
\cite[(8.2)(iii)]{hain-macp} for rational $n$-$K(\pi,1)$s
that if $X$ is a rational $n$-$K(\pi,1)$, then, for all $l$, the
natural
map
$$
W_lA_\C^t(X) \to W_lH^t(X;\C)
$$
is an isomorphism when $t\le n$ and injective when $t=n+1$.

The second step is more difficult. We will show that if $X$ is a
rational $n$-$K(\pi,1)$, then, for all $l$,
$$
W_lH^t(\Omega_\bw^\dot(X)_\F) \to W_lH^t(X,\R)
$$
is an isomorphism when $t\le n$ and injective when $t=n+1$.

First choose a base point $x\in X$. Let $\g$ be the complex form
of the Malcev Lie algebra associated to $(X,x)$. We shall view it
as a real Lie algebra with an almost complex structure $J$. Denote
its real form by $\g_\R^\blank$. We shall denote their (real)
continuous duals by $\g^\ast$ and $\g_\R^\ast$, respectively.

In \cite[(7.7)]{hain-macp}, a $\C$-linear map $\theta_x^\ast :
\Hom_\C(\g, \C) \to E^\dot_\C(X)$ is constructed and it is
established
that the image of $\theta^\ast_x$ is contained in $\Omega^\dot(X)$.
Define an $\R$ linear map
$$
\Theta_x : \g^\ast \to \Omega_\R^\dot(X)_\F
$$
as follows: each $\phi\in \g^\ast$ can be extended canonically to a
$\C$ linear map $\phihat : \g \to \C$. Define
$$
\Theta_x(\phi) = \Re \theta_x^\ast(\phihat).
$$
This induces an algebra homomorphism
$$
\Lambda^\dot_\R\g^\ast \to \Omega_\R^\dot(X)_\F.
$$
Since $\Theta$ clearly preserves the weight filtration, the induced
algebra homomorphism does too.

\begin{lemma}\label{first-qism}
Suppose that $X$ is a complex algebraic manifold with $q=0$. If $X$
is a rational $n$-$K(\pi,1)$, then, for all $l$,
$$
W_lH^t(\Lambda^\dot_\R \g^\ast) \to W_lH^t(\Omega^\dot_\R(X)_\F)
$$
is an isomorphism when $t\le n$ and injective when $t=n+1$.
\end{lemma}

\begin{pf}
Recall that $\g$ is viewed as a real Lie algebra with almost complex
structure $J$. Consequently,
$$
\g\otimes_\R \C = \g' \oplus \g''
$$
where $\g'$ and $\g''$, the $i$ and $-i$ eigenspaces of $J$,
respectively,
are commuting complex Lie subalgebras of $\g\otimes_\R\C$.. It
follows that
$$
\Lambda^\dot \g^\ast\otimes_\R \C \cong
\Lambda^\dot {\g'}^\ast \otimes \Lambda^\dot {\g''}^\ast.
$$
This is an isomorphism of $W_\dot$ filtered cochain complexes.

Similarly, there is an almost complex structure on the real vector
space $V$ with basis $\theta_j, \rho_j$, where $1\le j \le m$.  It is
defined by
$$
J : \theta_j \to \rho_j \text{ and } J : \rho_j \to -\theta_j.
$$
Define $V'$ and $V''$ to be the $i$ and $-i$ eigenspaces of $J$
acting
on $V\otimes \C$. Then $V'$ has basis
$$
df_j/f_j = \rho_j + i \theta_j, \quad j=1, \dots, m
$$
and $V'$ has basis their complex conjugates. It follows that
\begin{multline*}
\Omega_\R^\dot(X)_\F \cong \\
\Lambda_\C^\dot(df_j/f_j,j=1,\dots,m)/(a(u)+ ib(u)) \otimes
\Lambda_\C^\dot(d\fbar_j/\fbar_j,j=1,\dots,m)/(a(u)- ib(u)) \\
\cong A^\dot_\C(X) \otimes \overline{A^\dot_\C(X)}
\end{multline*}
Each of these algebras has the property that its degree $m$ part has
weight $2m$. Consequently, each of these isomorphisms is a
$W_\dot$ filtered algebra isomorphism.

The complexification of the map in the statement of the proposition
is the tensor product of the algebra homomorphism
$$
\Lambda^\dot {\g'}^\ast  \to A_\C^\dot(X)
$$
with its complex conjugate. It therefore suffices to prove that this
map is a $W_\dot$ filtered quasi-isomorphism.

Note that $\g'$ is just the complex form of the Malcev Lie algebra
associated to $\pi_1(X,x)$ and that the map above is the map
induced by the homomorphisms $\theta_x$ of
\cite[(7.7)]{hain-macp}.  This homomorphism induces a
homomorphism
\begin{equation}\label{qism}
\Lambda^\dot {\g'}^\ast \to
A^\dot_\C(X) \subseteq H^\dot(X;\C)
\end{equation}
which is the complexification of a morphism of mixed Hodge
structures. The result now follows as morphisms of mixed Hodge
structures are strict with respect to $W_\dot$ and since the map on
homology induced by (\ref{qism}) is an isomorphism in dimensions
$\le n$ and injective in dimension $n+1$ by the definition of a
rational $n$-$K(\pi,1)$.
\end{pf}

View $\Lambda^\dot_\R \g^\ast$ as the real left invariant
differential
forms on $G$, the complex form of the Malcev completion of
$\pi_1(X,x)$.
Here $G$ is viewed as a real proalgebraic group by restriction of
scalars. It follows from the fact that elements of $\bw(X)_\F$ are
represented by iterated integrals of elements of
$\Lambda^\dot_\R \g^\ast$ that the exterior derivative of each
element of $\bw(X)_\F$ is an element of
$\bw(X)_\F \otimes \Lambda^\dot_\R \g^\ast$
(cf. \cite[p.~436]{hain-macp}). Consequently,
$\bw(X)_\F \otimes \Lambda^\dot_\R \g^\ast$
is a subcomplex of $\lim_\to E^\dot_\R(G_s)$, the de~Rham complex of
$G$.

\begin{lemma}\label{second-qism}
Suppose that $X$ is a complex algebraic manifold with $q=0$. If $X$
is a rational $n$-$K(\pi,1)$, then the natural map
$$
\theta^\ast : \Omega_\bw^\dot(X)_\F \to E^\dot_\R(X)
$$
induces maps
$$
W_lH^m(\Omega_\bw^\dot(X)_\F) \to W_lH^m(X;\R)
$$
which are isomorphisms when $m\le n$ and injective when
$m=n+1$.
\end{lemma}

\begin{pf}
By considering the formal analogue of $\Omega_\bw^\dot(X)_\F$ for
$\Alb X$, it follows that one can put a differential on
$$
\bw(X)_\F \otimes \Lambda^\dot_\R \g^\ast
$$
such that $\Lambda^\dot_\R \g^\ast$ is a subcomplex and such that
the map to $\Omega_\bw^\dot(X)_\F$ induced by
$$
\Lambda^\dot_\R \g^\ast \to \Omega_\R^\dot(X)_\F
$$
is a $W_\dot$ filtered chain map. It follows from (\ref{first-qism})
that chain map
$$
\bw(X)_\F \otimes \Lambda^\dot_\R \g^\ast \to
\bw(X)_\F \otimes \Omega_\R^\dot(X)_\F = \Omega_\bw^\dot(X)_\F
$$
induces an isomorphism on $W_lH^k$ for all $l$ when $k\le n$ and
an injection on $W_lH^{n+1}$ for all $l$. (Filter each $Gr^W_l$ by
degree.)

To prove the result, we need only show that the complexification of
the composite
$$
\bw(X)_\F \otimes \Lambda^\dot_\R \g^\ast \otimes \C
\to E^\dot(\Xbar \log D)
$$
induces an isomorphism on $W_mH^t$ when $t\le n$ and an injection
on $W_mH^{n+1}$.

As in the proof of (\ref{first-qism}), we have the decomposition
$$
\g\otimes \C = \g' \oplus \g'',
$$
where $\g'$ and $\g''$ are commuting Lie algebras, and the $W_\dot$
filtered quasi-isomorphism
$$
\Lambda^\dot_\C \g \cong
\Lambda^\dot {\g'}^\ast \otimes \Lambda^\dot {\g''}^\ast
$$

Denote the Malcev group corresponding to $\g\otimes \C$ by $G_\C$,
and the commuting subgroups of $G_\C$ corresponding to $\g'$ and
$\g''$ by $G'$ and $G''$, respectively. Then
$$
G_\C = G' \times G''.
$$
Denote the complex points $G_\R(\C)$ of $G_\R(\C)$ by
$H$. Multiplication induces a continuous map
$$
H \times G'' \to G_\C.
$$
This is a continuous bijection.  To see this, let $H^s$, $G_\C^s$,
etc.
denote the $s$th terms of the lower central series of $H$, $G_\C$,
etc.
Note that the $s$th graded quotient of $\g_\C$ is the direct sum of
the $s$th graded quotients of the lower central series of $\h$ and
$\g''$. It follows that if $g \in G_\C$ is congruent to $hg''$ mod
$G_\C^s$, where $h\in H$ and $g''\in G''$, then there exit
$v\in H^s$, $u''\in {G'}^s$, unique mod $H^{s+1}$ and ${G''}^{s+1}$,
such that
$$
g(hg'')^{-1} = vu'' \text{ in } G^s_\C/G^{s+1}_\C.
$$
Since $G^s_\C/G^{s+1}_\C$ is central in $G_\C/G^{s+1}_\C$, it follows
that $g$ is congruent to $(hv)(g''u'')$ mod $G^{s+1}_\C$. The
assertion
follows by taking limits.

Since $\bw(X)_\F$ is the coordinate ring of the real proalgebraic
variety $G_\R\backslash G$, it follows that $\bw(X)_\F\otimes \C$ is
the coordinate ring of the complex proalgebraic variety $H\backslash
G_\C$. It follows immediately that the composite
$$
G'' \to G_\C \to H\backslash G_\C
$$
is an isomorphism of proalgebraic varieties. Consequently, there is a
natural algebra isomorphism
$$
\O(G'') \cong \bw(X)_\F \otimes \C.
$$
Since $G''$ is prounipotent, the exponential map induces
an isomorphism
$$
\C[{\g''}^\ast] \cong \O(G'').
$$

Assembling the pieces, we obtain an algebra isomorphism
\begin{equation}\label{filt_alg}
\bw(X)_\F \otimes \Lambda^\dot_\R \g^\ast \otimes \C
\cong \C[{\g''}^\ast] \otimes
\Lambda^\dot {\g'}^\ast \otimes \Lambda^\dot {\g''}^\ast
\end{equation}
The differential induced on the right hand side can be understood.
When the right hand side is quotiented out by the subcomplex
$\Lambda^\dot {\g'}^\ast$, the resulting complex is isomorphic to
the complex
$$
\C[{\g''}^\ast] \otimes  \Lambda^\dot {\g''}^\ast
$$
which is analogous to the complex in the proof of
\cite[(7.8)]{hain-macp}.  This complex is easily seen to by acyclic
by
the standard argument given there. It follows that the inclusion
\begin{equation}\label{alg_qism}
\Lambda^\dot {\g'}^\ast \hookrightarrow \C[{\g''}^\ast] \otimes
\Lambda^\dot {\g'}^\ast \otimes \Lambda^\dot {\g''}^\ast
\end{equation}
is a quasi-isomorphism.

It remains to show that this is a filtered quasi-isomorphism. First
observe that since $\g^\ast$ is the direct limit of complex parts of
mixed
Hodge structures (albeit, viewed as a real vector spaces), its weight
filtration has a canonical splittings by $J$ invariant subspaces. It
follows from the definitions that there are canonical splitting of
the
weight filtrations of ${\g'}^\ast$ and ${\g''}^\ast$ such that the
isomorphism
$$
\g^\ast \cong {\g'}^\ast \oplus {\g''}^\ast
$$
is an isomorphism of graded vector spaces. Consequently, there are
compatible
canonical splittings of the weight filtrations on
$$
\Lambda^\dot {\g'}^\ast,\quad \Lambda^\dot {\g''}^\ast, \quad
\Lambda^\dot_\C {\g}^\ast, \quad \C[{\g}^\ast], \quad \C[{\g'}^\ast],
\quad \C[{\g''}^\ast].
$$
Since $\bw(X)_\F = \R[{\g}^\ast]^\g$ and since the action
of $\g$ on $\R[{\g}^\ast]$ comes from a morphism of mixed Hodge
structures,
this action is compatible with the splittings of the weight
filtrations.
It follows that the weight filtration of $\bw(X)_\F$ has a
canonical splitting {\it that depends upon the choice of the base
point
$x$}. It follows that (\ref{filt_alg}) is an isomorphism of graded
algebras and that (\ref{alg_qism}) is a quasi-isomorphism of filtered
algebras, and therefore a $W_\dot$ filtered quasi-isomorphism.

To complete the proof, we have to show that the composite
$$
\Lambda^\dot {\g'}^\ast \hookrightarrow
\bw(X)_\F \otimes \Lambda^\dot_\R \g^\ast \otimes \C
\to E^\dot(\Xbar \log D)
$$
induces an isomorphism on $W_mH^t$ when $t\le n$ and injection on
$W_mH^{n+1}$.
Observe that if $\psi \in {\g'}^\ast$, then $\psi \in \g^\ast$ and
$\psi$
commutes with $J$. That is, $\psi \in \Hom_\C(\g,\C)$. It follows
that the
induced map
$$
\Lambda^\dot {\g'}^\ast\to E^\dot(\Xbar\log D)
$$
takes $\psi$ to $\theta_x^\ast(\psi)$. The assertion follows from the
fact that
$$
\theta^\ast_x : H^\dot(\g) \to H^\dot(X;\C)
$$
is a morphism of mixed Hodge structures \cite[(7.11)]{hain:cycles}
and that
$X$ is a rational $n$-$K(\pi,1)$.
\end{pf}

We are now ready to prove Theorem \ref{tech}. In the case where
$X_\dot$ is a single complex algebraic manifold with $q=0$, we have
the morphism
$$
\begin{matrix}
\to & F^{p}W_{2p}H^{t-1}(A^\dot_\C(X)) &
\to & \Hf^t(X,\R(m)) & \to & W_{2p}H^t(\Omega^\dot_\bw(X)_\F)
(m-1)  \cr
&\downarrow && \downarrow && \downarrow  \cr
\to & F^{p}W_{2p}H^{t-1}(X,\C) & \to & \Hd^t(X,\R(m))
& \to & W_{2p}H^t(X,\R(m-1))  \cr
\end{matrix}
$$
of long exact sequences. If $X$ is a rational $n$-$K(\pi,1)$, the
result
follows directly from Lemma \ref{second-qism}, the very first step
in the proof of Theorem \ref{tech}, and the 5-lemma.

In the simplicial case, the result follows by a similar argument:
Suppose that  $X_\dot$
is a simplicial complex algebraic manifold where each $X_m$ has
$q=0$. If each $X_m$ is a rational $(n-m)$-$K(\pi,1)$, then an
elementary
spectral sequence argument shows that the maps
$$
F^{p}W_{l}H^t(A^\dot_\C(X_\dot)) \to F^{p}W_{l}H^t(X_\dot;\C)
$$
and
$$
W_lH^t(\Omega^\dot_\bw(X_\dot)_\F) \to W_lH^t(X_\dot,\R)
$$
are isomorphisms for all $l$ and $p$ whenever $t\le n$. Theorem
\ref{tech} now follows from the 5-lemma as in the proof of the
result for a single space above.

\section{The 4-logarithm}

In this section, we prove the existence and uniqueness of a real
Grassmann 4-logarithm:

\begin{theorem}
There is a unique Grassmann 4-logarithm.
\end{theorem}

We will use the formal Bloch-Wigner
cohomology defined in Section \ref{comparison} as it behaves
better than $H^\dot_\bw$.

For all $m$, we have the following commutative diagram:
\begin{equation*}\label{ladder}
\begin{CD}
sW_{2m}H^{2m-1}(\Omega_\bw^\dot(G^m_\dot)_\F)
@>\alpha_{2m-1}>> sW_{2m}H^{2m-1}(G^m_\dot,\R) \\
@VVV @VVV \\
s\Hf^{2m}(G^m_\dot,\R(m)) @>>> s \Hd^{2m}(G^m_\dot,\R(m)) \\
@VVV @VVV \\
sF^mW_{2m}H^{2m}(A_\C^\dot(G^m_\dot)) @>>>
sF^mW_{2m}H^{2m}(\Omega^\dot(G^m_\dot))\\
@VVV @VVV\\
sW_{2m}H^{2m-1}(\Omega_\bw^\dot(G^m_\dot)_\F)
@>\alpha_{2m}>> sW_{2m}H^{2m}(G^m_\dot,\R) \\
\end{CD}
\end{equation*}
By (\ref{second-qism}), $\alpha_{2m-1}$ is an isomorphism if each
$G^m_n$  is a rational $(n-m-1)$-$K(\pi,1)$. If, in addition,
$$
sH^{m-n}(G^m_n,\R)
$$
vanishes when $n\ge 1$, then
$$
sH^{2m}(\Omega_\bw^\dot(G^m_\dot)_\F))
$$
is spanned by the class of $\vol_m$. Consequently,
$\alpha_{2m}$ is injective. Finally, observe that it follows from
\cite[(9.7)]{hain-macp} that the volume form
$\vol_m$ can be regarded as class in
$$
sF^mW_{2m}H^{2m}(A_\C^\dot(G^m_\dot))
$$

We now consider the case when $m=4$.
It follows from \cite[\S 8]{hain-macp} that $G^4_0$, $G^4_1$ are
rational $K(\pi,1)$s, and that $G^4_2$ is a rational 1-$K(\pi,1)$. It
follows from \cite[(8.2)]{hain-macp} that the cup products
$$
\Lambda^k H^1(G^4_n,\R) \to H^k(G^4_n,\R)
$$
are surjective when $k \le 4-n$, except possibly when $n=2$. But in
this case, the fibers of the face map $G^4_2 \to G^4_1$ are
hyperplane
complements with constant combinatorics.
It follows that $G^4_2 \to G^4_1$ is a fibration. The surjectivity of
the
cup product
$$
\Lambda^2 H^1(G^4_2) \to H^2(G^4_2)
$$
follows as the Leray-Serre spectral sequence of the map collapses at
$E^2$ for weight reasons. One can show (e.g., by computer or by
hand) that
$$
s\Lambda^k H^1(G^4_n,\R) = 0
$$
when $0 \le k < 4-n$ for all $n$, and when $k = 4-n$ for all $n> 0$.
By the discussion in the previous paragraph, $\alpha_4$ is injective
and
$$
sW_8H^7(\Omega_\bw^\dot(G^4_\dot)_\F) =
sW_8H^7(G^4_\dot,\R)=0.
$$
Thus, to prove the existence of the 4-logarithm, it suffices to show
that
$$
\vol_4\in sF^4W_8H^8(A_\C^\dot(G^4_\dot))
$$
has trivial image in
$sW_8H^7(\Omega_\bw^\dot(G^4_\dot)_\F)$. But this
is immediate as $\alpha_8$ is injective and $\vol_4$ has trivial
image in $H^8(G^4_\dot,\R)$.
The existence and uniqueness of the 4-logarithm follows.

\begin{remark}
It is likely that one can prove the existence of a canonical
5-logarithm using a similar argument.
\end{remark}

\end{document}